%% file: stability.tex
\documentclass[11pt,leqno]{article}
\usepackage{amsmath}
\usepackage{amssymb}
\usepackage{theorem}
\usepackage{curves}
\usepackage{epic}
\usepackage{rotating}
\usepackage{epsf}
\usepackage[scanall]{psfrag}
\usepackage{graphicx}

\numberwithin{equation}{section}

\newtheorem{theorem}{Theorem}[section]

\newtheorem{corollary}[theorem]{Corollary}
\newtheorem{lemma}{Lemma}[section]
\newtheorem{Claim}{Claim}
\theorembodyfont{\rmfamily}
\newtheorem{defn}{Definition}[section]

\newenvironment{proof}{\medskip\noindent{\bf Proof.}}{\medskip}

\newenvironment{mlem}{\medskip\noindent{\bf Main Lemma\
}\itshape}{\medskip}

\newcommand{\RR}{\mathbb{R}}
\newcommand{\po}{\partial\O}
\newcommand{\s}{\sigma}
\renewcommand{\H}{\mathcal H}
\newcommand\res{\hbox{ {\vrule height .22cm}{\leaders\hrule\hskip.2cm} } }
\newcommand{\cL}{\mathcal{L}}
\newcommand{\vecnQr}{\overrightarrow{n_{Q,r}}\,}
\newcommand{\vecnQs}{\overrightarrow{n_{Q,s}}\,}
\newcommand{\T}{\mathcal T}
\newcommand{\bs}{\backslash}

\renewcommand{\d}{\delta}

\renewcommand{\O}{\Omega}
\newcommand{\p}{\partial}
\newcommand{\loc}{\mathrm{loc}}
\newcommand{\NN}{\mathbb{N}}
\newcommand{\e}{\varepsilon}

\newcommand{\vecn}{\overrightarrow{n}\,}
\newcommand{\nint}{\int \kern-1.13em {\begin{turn}{-20}$\bigm/$%
\end{turn}}\!}
\newcommand\notint{{{\,\int \kern-1.01em \raise1pt\hbox{{\begin{turn}{-30}$/$%
\end{turn}\!\!}}}}}

\renewcommand{\SS}{\mathbb S}
\renewcommand{\a}{\alpha}
\renewcommand{\div}{\,\mathrm{div}\,}

\renewcommand{\o}{\omega}

\newcommand\qed{\hfill\vrule height8pt width6pt depth0pt}
\newcommand{\supp}{\mathrm{supp}}
\newcommand{\diam}{\mathrm{diam}}
\newcommand{\vmo}{\mathrm{VMO}}
\newcommand{\bmo}{\mathrm{BMO}}
\newcommand{\osc}{\mathrm{osc}\,}
\newcommand{\seg}{\mathrm{seg}\,}

\begin{document}

%\title{Stability of Lewis and Vogel's result}
%\author{D.\ Preiss and T.\ Toro}
%\date{}
%\maketitle

\input{stab-intro}

\input{stab-prel}
\input{stab-rough}

\input{stab-fine}
\input{stab-apps}

\input{stab-refs}

\end{document}

%% file: stab-intro.tex
\title{Stability of Lewis and Vogel's result}
\author{D.\ Preiss and T.\ Toro
\footnote{The second author was partially supported by NSF through DMS and 
by a Transitional Support Grant from ADVANCE at UW}}

\date{}
\maketitle

\section{Introduction}

Lewis and Vogel proved (see \cite{LV1}, \cite{LV2}) that a bounded
domain whose harmonic measure (with respect to a fixed point) is a
constant multiple of the surface measure to the boundary (i.e. a 
domain whose Poisson kernel is constant) is a ball, provided the surface
measure has at most Euclidean growth. In this paper we prove 
that this result is stable under small perturbations. Namely a bounded 
domain whose Poisson kernel is almost constant, and whose surface
measure to the boundary has at most Euclidean growth,
 is geometrically close 
to a ball. 

Both of these results can be viewed as free boundary regularity results 
for the Poisson kernel. An interesting feature is that 
regularity of the free boundary is proved without an a-priori assumption of 
flatness. In fact, our main theorem states that a domain whose Poisson kernel 
is almost constant has a locally flat boundary (see Theorem \ref{thm1}). Once 
the boundary is known to be locally flat the proof of regularity is standard.

Let $\O\subset\RR^{n+1}$ be a bounded domain and a set of locally finite
perimeter such that $0\in\O$ and $\H^n(\po)<\infty$. Let $\omega$ denote the
harmonic measure of $\O$ with pole at $0$. Let $\sigma$ denote the surface
measure of the boundary, i.e. $\s=\H^n\res\po$. Let $h=\frac{d\o}{d\s}$
denote the Poisson kernel of $\O$ with pole at $0$. First we state Lewis
and Vogel's result. Then we state one of our results which emphasizes
the stability of their result.

\begin{theorem}\label{thm1A}{\bf{\cite{LV1}}}
Assume that $\O\subset\RR^{n+1}$ satisfies 
%\begin{itemize}
%\item{}
\begin{equation}\label{eqn1A}
\sup_{0<r<1}\sup_{Q\in\po} \frac{\H^n(B(Q,r)\cap\po)}{r^n}<\infty,
\end{equation}
%\item{}
\begin{equation}\label{eqn2A}
\omega=\H^n\res \po.
\end{equation}
%\end{itemize}
Then $\O$ is a ball of center $0$ and radius $R>0$ such that 
$\H^n(\p B(0,R))=1$.
\end{theorem}

\begin{theorem}\label{thm2A}
Assume that $\O\subset\RR^{n+1}$ satisfies 
%\begin{itemize}
%\item{}
\begin{equation}\label{eqn3A}
\sup_{0<r<1}\sup_{Q\in\po} \frac{\H^n(B(Q,r)\cap\po)}{r^n}<\infty,
\end{equation}
%\item{}
\begin{equation}\label{eqn4A}
\frac{d\omega}{d\H^n}=h \qquad {\rm and} \qquad \sup_{\po}|\log h|<\e,
\end{equation}
for some $\e>0$ small enough,
%\item{} $\log h\in C^{\infty}$.
%\end{itemize}
Then $\O$ is a ``smooth'' deformation of $B(0,R)$
and $D[B(0,R),\O]<4\e$. Here $\H^n(\p B(0,R))=1$ and
$D$ denotes the Hausdorff distance.
\end{theorem}

The paper is organized as follows: in section 2 we introduce 
some definitions and state  
the main theorem precisely. In section 3 we 
prove that the gradient of the Green function near the boundary 
is controlled by the Poisson kernel. This is a consequence of the fact that 
the gradient of the Green function is a subharmonic function on a bounded 
domain and therefore the values near the boundary are controlled
by the boundary values. Recall that the Poisson kernel is basically the 
derivative of the Green function at the boundary.
As a consequence we show that if $\O$ satisfies
(\ref{eqn3A}) and (\ref{eqn4A}) then $D[B(0,R),\O]<4\e$. In section 
4 we introduce 
a local notion of flatness which involves the geometry of the boundary 
at a point and the behavior of $G$ and $\log h$ near that point. This
allows us to show that $\po$ is locally flat. In section 5 we present 
some applications of Theorem \ref{thm1}.

%% file: stab-prel.tex
\section{Preliminaries}

In this section we introduce the definitions needed to state our main results.
The main theorem appears at the end of the section and it is proved in section 4. We always assume that
$n\ge 2$.

\begin{defn}\label{defn1.1}
Let $\Sigma\subset\RR^{n+1}$ be a locally compact set, and let $\d>0$. We
say that $\Sigma$ is \emph{$\d$-Reifenberg flat} if for each compact set
$K\subset\RR^{n+1}$, there exists $R_K>0$ such that for every $Q\in
K\cap\Sigma$ and every $R\in (0, R_K]$ there exists an $n$-dimensional plane
$L(Q,r)$ containing $Q$ such that
\begin{equation}\label{eqn1.1}
\frac{1}{r}D[\Sigma\cap B(Q,r),L(Q,r)\cap B(Q,r)]\le\d.
\end{equation}
Here $B(Q,r)$ denotes the $(n+1)$-dimensional ball of radius $r$ and center
$Q$, and $D$ denotes the Hausdorff distance.
\end{defn}

Recall that for $A, B\subset\RR^{n+1}$,
\[
D[A,B]=\sup\{d(a,B):a\in A\}+\sup\{d(b,A):b\in B\}.
\]
Note that the previous definition is only significant for $\d>0$ small. We
denote by
\begin{equation}\label{eqn1.2}
\theta(Q,r)=\inf_{L}\left\{\frac{1}{r} D[\Sigma\cap B(Q,r),L\cap B
(Q,r)]\right\},
\end{equation}
where the infimum is taken over all $n$-planes containing $Q$.

\begin{defn}\label{defn1.4}
Let $\O\subset\RR^{n+1}$ be a set of locally finite perimeter (see \cite{EG}),
$\po$ is said to be Ahlfors regular if the surface measure to the boundary,
i.e., the restriction of the $n$-dimensional Hausdorff measure to $\po$,
$\s=\H^n\res \po$, is Ahlfors regular. That is there exists a constant
$C>1$
so that for $Q\in\po$ and $r\in(0, \rm{diam}\O)$
\begin{equation}\label{eqn1.3}
C^{-1}r^n\le\s(B(Q,r))\le Cr^n.
\end{equation}
\end{defn}

\begin{defn}\label{defn1.5}
Let $\O\subset\RR^{n+1}$ be a bounded set. We say that $\O$ 
has the \emph{separation
property} if there exists $R>0$
such that for $Q\in\po$ and $r\in(0,R]$ there exists an
$n$-dimensional plane $\cL(Q,r)$ containing $Q$ and a choice of unit normal
vector to $\cL(Q,r),\vecnQr$ satisfying
\begin{equation}\label{eqn1.4}
\T^+(Q,r)=\left\{X=(x,t)=x+t\vecnQr\in
B(Q,r):x\in\cL(Q,r),t>\frac{1}{4}r\right\}\subset\O,
\end{equation}
and
\begin{equation}\label{eqn1.5}
\T^-(Q,r) = \left\{X=(x,t)=x+t\vecnQr\in
B(Q,r):x\in\cL(Q,r),t<-\frac{1}{4}r\right\}\subset\O^c.
\end{equation}
\end{defn}

The notation $(x,t)=x+t\vecnQr$ is used to denote a point in $\RR^{n+1}$.
The first component, $x$, of the pair belongs to an $n$-dimensional affine
space whose unit normal vector is $\vecnQr$. The second component $t$
belongs to $\RR$. From the context it will always be clear what affine
hyperplane $x$ belongs to, and what the orientation of the unit normal
vector is.

\begin{defn}\label{defn1.6}
Let $\d\in(0,\d_n)$, where $\d_n$ is chosen appropriately (see note below)
and let $\O\subset\RR^{n+1}$. We say that $\O$ is a \emph{$\d$-Reifenberg flat
domain} or a \emph{Reifenberg flat domain} if $\O$ has the separation
property and $\po$ is $\d$-Reifenberg flat.  
\end{defn}

When we consider
$\d$-Reifenberg flat domains in $\RR^{n+1}$ we assume that $\d_n>0$ is small
enough, in order to ensure that we are working on NTA domains (see
definition in Appendix A, see also \cite{JK} and \cite[Theorem
3.1]{KT2}).

\begin{defn}\label{defn1.8}
A set of locally finite perimeter $\O\subset\RR^{n+1}$ is
said to be a \emph{chord arc domain}, if $\O$ is an NTA domain whose
boundary is Ahlfors regular.
\end{defn}

\begin{defn}\label{defn1.9}
Let $\d\in(0,\d_n)$. A set of locally finite perimeter $\O\subset\RR^{n+1}$
is said to be a \emph{$\d$-Reifenberg flat chord arc domain}, if $\O$ is a
$\d$-Reifenberg flat domain whose boundary is Ahlfors regular.
\end{defn}

\begin{defn}\label{defn1.10}
Let $\d\in(0,\d_n)$. A bounded set of locally finite perimeter $\O$
is said to be a \emph{$\d$-chord arc domain} or a \emph{chord arc domain
with small constant} if $\O$ is a $\d$-Reifenberg flat domain, $\po$ is
Ahlfors regular and there exists
$R>0$ so that
\begin{equation}\label{eqn1.11}
\sup_{Q\in\po\cap K}\|\vecn\|_\ast(Q,R)<\d.
\end{equation}
Here $\vecn$ denotes the unit normal vector to the boundary,
\begin{equation}\label{eqn1.12}
\|\vecn\|_\ast(Q,R)=\sup_{0<s<R}
(\notint_{B(Q,s)}|\vecn-\vecnQs|^2d\s)^{\frac{1}{2}}
\end{equation}
and $\vecnQs=\notint_{B(Q,s)}\vecn d\s$.
\end{defn}

\begin{defn}\label{defn1.12}
Let $\O\subset\RR^{n+1}$ be a chord arc domain. Let $f\in L^2_{\loc}(d\s)$,
we say that $f\in \bmo(\po)$ if
\begin{equation}\label{eqn1.14}
\|f\|_\ast=\sup_{r>0} \sup_{Q\in\po}
(\notint_{B(Q,r)}|f-f_{Q,r}|^2d\s)^{\frac{1}{2}}<\infty.
\end{equation}
Here $f_{Q,r}=\notint_{B(Q,r)}fd\s$, and $\s=\H^n\res\po$.
\end{defn}

\begin{defn}\label{defn1.13}
Let $\O\subset\RR^{n+1}$ be a chord arc domain. We denote by $\vmo(\po)$ the
closure in $\bmo(\po)$ of the set of uniformly continuous bounded functions
defined on $\po$.
\end{defn}

From now on we assume that 
$\O\subset\RR^{n+1}$ is a bounded domain and a set of locally finite
perimeter such that $0\in\O$ and $\H^n(\po)<\infty$. Let $\omega$ denote the
harmonic measure of $\O$ with pole at $0$. Let $\sigma$ denote the surface
measure of the boundary. Let $h=\frac{d\o}{d\s}$
denote the Poisson kernel of $\O$ with pole at $0$.

\begin{theorem}\label{thm1}
Assume that $\O\subset\RR^{n+1}$ satisfies 
\begin{equation}\label{eqn1}
\sup_{0<r<1}\sup_{Q\in\po} \frac{\H^n(B(Q,r)\cap\po)}{r^n}<\infty.
\end{equation}
Then given $\sigma>0$ small enough there exists $\e>0$ such that if
\begin{equation}\label{eqn2}
\sup_{\po}|\log h|<\e
\end{equation}
then $\po$ is $\s$-Reifenberg flat.
\end{theorem}

%% file: stab-rough.tex
\section{Rough geometric properties}

The Main Lemma below provides a crucial estimate of the 
gradient of the Green function near the boundary 
in terms of the Poisson kernel. It allows us to deduce that 
under the hypothesis of Theorem \ref{thm1}, $\po$ is contained 
in a very thin annular region.

\begin{mlem}
Let $\O\subset\RR^{n+1}$, $0\in\O$. Let $G$ denote the Green function of
$\O$ with pole $0$ and let $h$ be the corresponding Poisson kernel. Assume
that
\begin{equation}\label{eqn15}
\sup_{0<r<1}\sup_{Q\in\po} \frac{\H^n(B(Q,r)\cap\po)}{r^n}<\infty
\end{equation}
and
\begin{equation}\label{eqn16}
\sup_{\po}|\log h|<\e
\end{equation}
for some $\e\in(0,1)$. Then
\begin{equation}\label{eqn17}
\limsup_{X\to P} |\nabla G(X)|\le e^\e\qquad\forall\;P\in\po.
\end{equation}
\end{mlem}

Let 
\begin{equation}\label{eqn18}
K_0=\sup_{0<r<1} \sup_{Q\in\po} \frac{\H^n(B(Q,r)\cap\po)}{r^n}<\infty.
\end{equation}

\begin{lemma}\label{lem3}
Under the assumptions above, let $R>0$ be such that $B(0,R)\subset\O$ and
$\p B(0,R)\cap\po\ne\emptyset$. Then
\begin{equation}\label{eqn19}
|\nabla G(X)|\le C_nK_0\qquad\forall\,X\in\O\backslash
B\left(0,\frac{R}{2}\right).
\end{equation}
\end{lemma}

\begin{proof}
Apply the Riesz decomposition theorem for subharmonic functions to $G$ (see
\cite[Theorem 6.18]{H}). Let $Q\in\po$ be such that $0\not\in
B(Q,r)$
\begin{eqnarray}\label{eqn20}
G(Q) & = & \notint_{\lower.1in\hbox{$\scriptstyle{\p B(Q,r)}$}}G(Z)d\sigma(Z) \\
&& \quad
-\frac{1}{(n-1)(n+1)\o_{n+1}}\int_{B(Q,r)\cap\po}\left(\frac{1}{|Z-Q|^{n-1}}
- \frac{1}{r^{n-1}}\right)d\s(Z).\nonumber
\end{eqnarray}
Using Fubini and the fact that $G(Q)=0$ (\ref{eqn20}) yields
\begin{equation}\label{eqn21}
\notint_{\lower.1in\hbox{$\scriptstyle{\p B(Q,r)}$}} G(Z)d\s(Z) = \frac{1}{(n+1)\o_{n+1}}\int^r_0
\frac{\o(B(Q,t))}{t^n}dt.
\end{equation}
Note that (\ref{eqn2}) and (\ref{eqn18}) imply that for $t<1$,
\begin{equation}\label{eqn22}
\o(B(Q,t)) \le e^\e\H^n(B(Q,t)\cap\po) \le  e^\e K_0t^n. 
\end{equation}
Combining (\ref{eqn21}) and (\ref{eqn22}) we have that for $\e<1$
\begin{equation}\label{eqn23}
\notint_{\lower.1in\hbox{$\scriptstyle{\p B(Q,r)}$}}G(Z)d\s(Z)\le C_nK_0r
\end{equation}
whenever $Q\in\po$ and $0\not\in B(Q,r)$.

Let $X\in\O\backslash B\left(0,\frac{R}{4}\right)$, there exists $Q\in\po$
such that $d(X)=r=|X-Q|$ where $d(X)$ denotes the distance from $X$ to
$\po$. If $r<\frac{R}{4}$ then $0\not\in B(Q,4r)$, and the representation
formula for subharmonic functions implies
\begin{equation}\label{eqn24}
G(X)\le \frac{(2r)^2-|X-Q|^2}{(n+1)\o_{n+1}(2r)} \int_{\p B(Q,2r)}
\frac{G(Z)}{|Z-X|^{n+1}}d\s(Z).
\end{equation}
Since $|Z-X|\ge r$ for $X\in\p B(Q,r)$, (\ref{eqn23}) and (\ref{eqn24})
yield
\begin{equation}\label{eqn25}
G(X) \le \frac{3}{2} \notint_{\lower.1in\hbox{$\scriptstyle\p B(Q,2r)$}} G(Z)d\s(Z)\le C_nK_0r = C_nK_0d(X).
\end{equation}

If $d(X)=r>\frac{R}{4}$, $X\in B\left(Q, \frac{7}{8}R\right)$ and $0\not\in
B\left(Q, \frac{7}{8}R\right)$. A similar argument to the one sketched above
proves that
\begin{equation}\label{eqn26}
G(X)\le C_n\notint_{\lower.1in\hbox{$\scriptstyle{\p
B\left(Q,\frac{7}{8}r\right)}$}} G(Z)d\s(Z) \le C_nK_0R
\le C_nK_0d(X).
\end{equation}
Thus we have shown that for $X\in\O\backslash B\left(0,\frac{R}{4}\right)$
\begin{equation}\label{eqn27}
G(X)\le C_nK_0d(X).
\end{equation}
Standard estimates for harmonic functions on $\O\backslash B\left(0,
\frac{R}{2}\right)$ ensure that
\begin{equation}\label{eqn28}
|\nabla G(X)|\le C_n \frac{G(X)}{d(X)} = C_nK_0
\end{equation}
\qed
\end{proof}

The proof of the Main Lemma is a slight variation of the proof that appears
in \cite{LV1}. We sketch the proof and try
to indicate as we go along what the ideas behind the calculations are. For
further details we refer the reader to \cite{LV1} and \cite{LV2}.
\medskip

\noindent{\bf Proof of Main Lemma:}  Let
$M=\limsup_{X\to\po}|\nabla G(X)|$. Assume that $M>e^\e$.
Let $\d\in(0,10^{-10})$ and let $X_0\in\O\bs B\left(0,\frac{3R}{4}\right)$
be such that
\begin{equation}\label{eqn29}
|\nabla G(X_0)|\ge M-\d.
\end{equation}
Let $W(X)=\max\{|\nabla G(X)|-(M-2\d); 0\}$, observe that $W(X_0)\ge \d$,
and that $W$ is subharmonic in $\O\bs B\left(0, \frac{R}{2}\right)$. Let
$G_0$ be the Green's function of $\O$ with pole at $X_0$. By Sard's theorem
we can choose $t>0$ such that $|\nabla G_0(X)|\ne 0$ on $\{X:G_0(X)=t\}$.
Green's second identity, the fact that $W$ is subharmonic on $\O\bs
B\left(0, \frac{R}{2}\right)$, the maximum principle applied to $G$ and
$G_0$ on $\O\bs B\left(0, \frac{R}{2}\right)$ and $\O\bs B\left(X_0,
\frac{d_0}{2}\right)$ respectively, where $d_0=d(X_0)$  and (\ref{eqn19})
yield
\begin{equation}\label{eqn30}
\frac{1}{6}\le \int_{\{|\nabla G|>M-2\d\}\cap\{G_0=t\}} \frac{\p
G_0}{\p\nu} (u)d\H^n(Y)
\end{equation}
provided $X_0$ is close enough to $\po$, and $t$ is chosen small enough so
that $|\nabla G|<M+\d$ on $\{X:G_0(X)=t\}$. Let $E(t)=\{X:|\nabla
G(X)|>M-2\d\}\cap\{X:G_0(X)=t\}$.

First one shows that $E(t)$ is a ``large'' set at ``distance'' comparable to
$t$ from $\po$. More precisely for $X\in E(t)$,
\begin{equation}\label{eqn31}
C_1d(X)\le t\le C_2d(X)
\end{equation}
where
$C_i=C(n,K_0,R,X_0)$ for $i=1,2$. Furthermore for $t$ small enough there exist balls
$\{B(X_i, d(X_i)\}$ with $X_i=X_i(t)\in E(t)$ such that
\begin{gather}
E(t)\subset\bigcup_i B\left(X_i, \frac{d(X_i)}{4}\right) \label{eqn32} \\
B\left(X_i, \frac{d(X_i)}{100}\right)\cap B\left(X_j,
\frac{d(X_j)}{100}\right) = \emptyset\quad\mbox{for }i\ne j \label{eqn33} \\
\sum_id(X_i)^n\ge C_3^{-1}, \label{eqn34}
\end{gather}
where $C_3=C(X_0,K_0, n, R)$. Note that each $B(X_i, d(X_i))$ is tangent to
$\po$.

Let $\gamma>0$ be a small positive constant. Since $\O$ is a set of locally
finite perimeter, Egoroff's theorem ensures that there exits $r_\gamma>0$ so
that
\begin{equation}\label{eqn35}
\frac{\H^n(\po\cap B(Z,r))}{\o_nr^n} <1+\gamma\ \mbox{ for }\ 0<r<r_\gamma
\end{equation}
whenever $Z\in\po\bs \Lambda$ and $\H^n(\Lambda)<\gamma^{100n}$. Choosing
$t\ll r_\gamma$ (\ref{eqn32}), (\ref{eqn33}), (\ref{eqn34}) and Lemma 3 in
\cite{LV1} guarantee that there exists $Y\in E(t)$ so that
\begin{equation}\label{eqn36}
|\nabla G(X)-\nabla G(Y)|\le\gamma\quad\forall\;X\in B(Y; (1-\gamma)d(Y))
\end{equation}
and if $\widehat Z\in\po\cap \p B(Y, d(Y))$ then there exists $Z\in\po$ such
that $|Z-\widehat Z|<\gamma t$ and $Z$ satisfies (\ref{eqn35})

For $0<r<r_0$ (\ref{eqn19}), (\ref{eqn21}), (\ref{eqn22}), (\ref{eqn36}) and
the fact that $t\sim d(Y)$ yield
\begin{eqnarray}\label{eqn37}
\notint_{\lower.1in\hbox{$\scriptstyle{\p B(\widehat Z,r)}$}}Gd\sigma & \le &
\notint_{\lower.1in\hbox{$\scriptstyle{\p B(Z,r)}$}} Gd\s +
C_nK_0\gamma t \\
& \le & \frac{1}{(n+1)\o_{n+1}} \int^r_0 \frac{\o(B(Z,s))}{s^n} ds +
C_nK_0\gamma d(Y) \nonumber \\
& \le & \frac{e^\e}{(n+1)\o_{n+1}} \int^r_0 \frac{\H^n(B(Z,s)\cap\po)}{s^n}
ds+C_nK_0\gamma d(Y) \nonumber \\
& \le & \frac{e^\e\o_n}{(n+1)\o_{n+1}} (1+\gamma)r + C_nK_0\gamma d(Y)
\nonumber
\end{eqnarray}
Assume $\widehat Z=Y-d(Y)e$, from (\ref{eqn19}) and (\ref{eqn35}) we deduce
for $X\in B(Y, d(Y))$
\begin{equation}\label{eqn38}
|G(X)-G(Y)-\langle\nabla G(Y); X-Y\rangle|\le C_nK_0\gamma d(Y).
\end{equation}
For $X=\widehat Z$ we have
\begin{equation}\label{eqn39}
|G(Y)-\langle\nabla G(Y); d(Y)e\rangle|\le C_nK_0\gamma d(Y).
\end{equation}
Combining (\ref{eqn38}) and (\ref{eqn39}) and using the fact that $G\ge 0$
we obtain for $X\in B(Y,d(Y))$
\begin{equation}\label{eqn40}
-\langle\nabla G(Y),e\rangle d(Y) - \langle\nabla G(Y); X-Y\rangle \le
2C_nK_0\gamma d(Y).
\end{equation}
Since $Y\in E(t)$, $|\nabla G(Y)|\ne 0$, letting $X$ tend to $-d(Y)\nabla
G(Y)/|\nabla G(Y)|$ we obtain
\begin{equation}\label{eqn41}
0\le|\nabla G(Y)|-\langle \nabla G(Y),e\rangle\le 2C_nK_0\gamma.
\end{equation}
Combining (\ref{eqn38}), (\ref{eqn39}) and (\ref{eqn40}) we find that
\begin{equation}\label{eqn42}
|G(X)-\langle\nabla G(Y); e\rangle(\langle X-Y,e\rangle+d(Y))|\le
C_nK_0\gamma d(Y)
\end{equation}
for $X\in B(Y, d(Y))$. Let $r=\gamma^{1/2}d(Y)$. Note that 
\begin{equation}\label{eqn43}
\H^n(\{X:\langle X-Y,e\rangle +d(Y)\ge 0\}\bs B(Y,d(Y))\cap \p B(\widehat
Z,r))\le C\gamma^{1/2}r^n
\end{equation}
and on this set
\begin{equation}\label{eqn44}
\langle X-Y;e\rangle + d(Y)\le C\gamma^{1/2}r.
\end{equation}
From (\ref{eqn42}), (\ref{eqn43}), (\ref{eqn44}) and (\ref{eqn19}) we have
\begin{eqnarray}\label{eqn45}
&&\kern-.5in\int_{\p B(\widehat Z,r)} G(X)d\s(X) \\[3mm]
& \ge & \int_{\p B(\widehat
Z,r)\cap\{X:\langle X-Y,e\rangle+d(Y)\ge 0\}}G(X)d\s(X)\nonumber \\[3mm]
& \ge & \int_{\p B(\widehat Z,r)\cap\{X:\langle X-Y,e\rangle +d(Y)\ge 0\}\cap
B(Y,d(Y))}G(X)d\s(X) \nonumber \\[3mm]
& \ge & \langle\nabla G(Y),e\rangle\int_{\p B(\widehat Z,r)\cap\{X:\langle
X-Y,e\rangle+d(Y)\ge 0\}\cap B(Y, d(Y))} (\langle X-Y,e\rangle
 +d(Y))d\s(X) \nonumber \\[3mm]
& \ge & \langle\nabla G(Y),e\rangle\int_{\p B(\widehat Z,r)\cap\{X:\langle
X-Y,e\rangle+d(Y)\ge 0\}} (\langle X-Y,e\rangle + d(Y))d\s(X) \nonumber \\[3mm]
&& \qquad\qquad\qquad\qquad\qquad\qquad\qquad\qquad\qquad
 -C|\langle\nabla G(Y),e\rangle|\gamma^{1/2}r^{n+1} 
\nonumber \\[3mm]
& \ge & \langle\nabla G(Y),e\rangle \int_{\p B(\widehat Z,r)\cap\{X:\langle
X-Y,e\rangle + d(Y)\ge 0\}} (\langle X-Y;e\rangle+d(Y))d\s(X)\nonumber \\
&&\qquad\qquad\qquad\qquad\qquad\qquad\qquad\qquad\qquad
 -C_nK_0\gamma^{1/2}r^{n+1}.\nonumber
\end{eqnarray}

Note that
\begin{equation}\label{eqn46}
\int_{\p B(\widehat Z,r)\cap\{X:\langle X-Y;e\rangle+d(Y)\ge 0\}}(\langle
X-Y,e\rangle+d(Y))d\s(X) = \int_{\p B(0,r)\cap\{X:x_{n+1}\ge
0\}}x_{n+1}d\s(X).
\end{equation}
The representation formula for subharmonic functions applied to
$V(X)=\max\{x_{n+1},0\}$ yields
\begin{equation}\label{eqn47}
\int_{\p B(0,r)\cap\{X:x_{n+1}\ge 0\}}x_{n+1}d\s(X) = r^n\int^r_0
\frac{\o_ns^n}{s^n}ds = \o_nr^{n+1}.
\end{equation}

Combining (\ref{eqn45}), (\ref{eqn46}) and (\ref{eqn47}) we have 
\begin{equation}\label{eqn48}
\int_{\p B(\widehat Z,r)} G(X)d\s(X) \ge \langle\nabla G(Y),e\rangle
\o_nr^{n+1} - C_nK_0\gamma^{1/2}r^{n+1}.
\end{equation}
From (\ref{eqn37}) and (\ref{eqn48}) we deduce
\begin{equation}\label{eqn49}
\langle\nabla G(Y),e\rangle \frac{\o_n}{(n+1)\o_{n+1}}r-C_nK_0\gamma^{1/2}r
\le \frac{\o_n}{(n+1)\o_n} e^\e(1+\gamma)r + C_nK_0\gamma^{1/2}r
\end{equation}
thus
\begin{equation}\label{eqn50}
\langle\nabla G(Y), e\rangle\le e^\e(1+\gamma)+C_nK_0\gamma^{1/2}.
\end{equation}
Using the fact that $Y\in E(t)$, (\ref{eqn41}) and (\ref{eqn50}) we conclude
that
\begin{equation}\label{eqn51}
M-2\d\le |\nabla G(Y)|  \le  \langle\nabla G(Y),e\rangle +
2C_nK_0\gamma 
 \le  e^\e(1+\gamma)+C_nK_0\gamma^{1/2}. 
\end{equation}
Since $\gamma>0$ is arbitrary we conclude
from (\ref{eqn51}) that $M-2\d\le e^\e$. Letting $\d$ tend to $0$ we get
that $M\le e^\e$, which contradicts our initial assumption that $M>e^\e$.
This remark finishes the proof of the main lemma.\hfill\qed

Let
\begin{equation}\label{eqn3}
0<R_1=\sup\{r:B(0,r)\subset\O\}<\infty
\end{equation}
\begin{equation}\label{eqn4}
0<R_2=\inf\{r:\O\subset B(0,r)\}<\infty.
\end{equation}
To estimate $R_1$, let $P_1=\po\cap\p B(0,R_1)$. Let $G_1$ be the
Green's function of $B(0,R_1)$ with pole $0$, let $G$ be the Green's
function of $\O$ with pole $0$. By the maximum principle for $X\in
B(0,R_1)\backslash\{0\}$
\begin{equation}\label{eqn5}
G_1(X)\le G(X).
\end{equation}
In fact if $F(X)$ denotes the fundamental solution for the Laplacian in
$\RR^{n+1}$ with pole at the origin then $G=F-u$ and $G_1=F-u_1$ where
$\Delta u=0$ in $\O$ with $u=F$ on $\po$ and $\Delta u_1=0$ in $B(0,R_1)$
with $u_1=F$ on $\p B(0,R_1)$. Since $G\ge 0$ then $u\le F$ in $\O$, and
hence $u\le u_1$ on $\partial B(0,R_1)$ (because $B(0,R_1)\subset\O$). By
the maximum principle $u\le u_1$ in $B(0,R_1)$ which justifies
(\ref{eqn5}). Letting $X=tP_1$ with $t\to 1$ (\ref{eqn5})
yields
\begin{equation}\label{eqn6}
\liminf_{t\to 1} \frac{G_1(tP_1)}{t} \le \liminf_{t\to 1} \frac{G(tP_1)}{t}.
\end{equation}
Thus by (\ref{eqn2}) and the Main Lemma we have that
\begin{equation}\label{eqn7}
\frac{1}{\H^n(\p B(0,R_1))}=|\nabla G_1(P_1)|\le e^\e.
\end{equation}
If $\s_n=\H^n(\p B(0,1))$ then (\ref{eqn7}) implies
\begin{equation}\label{eqn8}
(e^{-\e}\s_n^{-1})^{1/n}\le R_1.
\end{equation}
To estimate $R_2$ let $P_2\in\po$ be such that $|P_2|=\max\{|Q|:Q\in\po\}$.
Let $G_2$ denote the Green's function of $B(0,R_2)$ with pole at $0$. A
similar argument to the one above shows that for $X\in\O\backslash\{0\}$
\begin{equation}\label{eqn9}
G(X)\le G_2(X).
\end{equation}
Note that for $P_2$ there exists a ball $B\subset\O^c$ such that
$P_2\in\po\cap\p B$.

\begin{lemma}\label{lemma1}
Let $\O$, $G$ and $h$ be as above. Let $P\subset\po$ and assume that there
exists a ball $B\subset\O^c=\{G=0\}$ so that $P\in\po\cap\p B$ then
\begin{equation}\label{eqntilde1}
\limsup_{\mathop{X\to P}\limits_{X\in\O}} \frac{G(X)}{d(X,B)}\ge e^{-\e}
\end{equation}
\end{lemma}

\begin{proof}
Let $l=\mathop{\limsup}\limits_{\mathop{X\to P}\limits_{X\in\O}}\frac{G(X)}{d(X,B)}$. There
exists a sequence $\{Y_k\}_{k\ge 1}\subset\O$, such that $Y_k\to P$ and
$\frac{G(Y_k)}{d(Y_k,B)}\mathop{\longrightarrow}\limits_{k\to\infty} l$. Let
$d_k=d(Y_k,B)$. There exists $X_k\in\p B$ so that $|Y_k-X_k|=d_k$. Consider
$G_k(X)=\frac{G(d_kX+X_k)}{d_k}$ for $X\in B(0,2)$ and
$Z_k=\frac{Y_k-X_k}{d_k}$. 
Without loss of generality we may assume that $Z_k\to e$ as $k\to\infty$,
$|e|=1$, and $G_k\mathop{\longrightarrow}\limits_{k\to\infty}G_\infty$ in
$C^{0,\beta}_{\loc}(\RR^{n+1})$, $\nabla G_k
\mathop{\rightharpoonup}\limits_{k\to\infty}^{\ast} \nabla G_\infty$ weak
star in $L^\infty_{\loc}(\RR^{n+1})$, weakly in $L^2_{\loc}(\RR^{n+1})$;
$\frac{1}{d_k}(\po-X_k)=\p\{G_k>0\}
\mathop{\longrightarrow}\limits_{k\to\infty}\p\{G_\infty>0\}$ in the
Hausdorff distance sense uniformly on compact sets, and $\chi_{\{G_k>0\}}\to
\chi_{\{G_\infty>0\}}$ in $L^1_{\loc}(\RR^{n+1})$. Note that
$G_k(Z_k)=\frac{G(Y_k)}{d_k}$ thus $G_k(Z_k)\to l$ as $k\to\infty$. On the
other hand since $G_k$ converges uniformly to $G_\infty$ in $B(0,2)$, we
conclude that $G_\infty(e)=l$. In order to prove the lemma we need to get a
better understanding of $G_\infty$ and $\O_\infty=\{G_\infty>0\}$. Our goal
is to show that $\O_\infty$ is a half-space and $G_\infty$ is linear. 
Let $r$ be
the radius of $B$.
Let $\a_k=d(\p B(X_k,d_k)\cap\p B; L)$, where $L$ is the tangent plane to
$B$ through $X_k$. An easy computation shows that $\a_k=2\frac{d^2_k}{r}$.
Note that for $P_k\in B(X_k,d_k)\cap\{\langle P-X_k,
\frac{Y_k-X_k}{d_k}\rangle <-\a_k\}\subset B$ if $Q_k=\frac{P_k-X_k}{d_k}$,
then $Q_k\in B(0,2)\cap\{\langle X, Z_k\rangle < -\frac{d_k}{r}\}$ and
$G_k(Q_k)\le 0$. Passing to the limit as $k$ tends to infinity we conclude
that if $Y\in B(0,2)\cap\{\langle Y,e\rangle\le 0\}$ then $G_\infty(Y)=0$.
Let $Y\in B(0,2)\cap\{\langle Y, Z_k\rangle >0\}$, then either
$d_kY+X_k\in\O^c$ and $G_k(Y)=0$ or $d_kY+X_k\in \O$ and given $\e>0$ there
exists $k_0\in\NN$ such that for $k\ge k_0$
\begin{equation}\label{eqntilde2}
\frac{G(d_kY+X_k)}{d(d_kY+X_k,B)}\le l+\e
\end{equation}
and
\begin{eqnarray}\label{eqntilde3}
G(d_kY+X_k) & \le & (l+\e)d(d_kY+X_k,B) \\
& \le & (l+\e)\left\{\langle d_kY;
\frac{Y_k-X_k}{d_k}\rangle+2\frac{d^2_k}{r}\right\} \nonumber \\
& \le & (l+\e)d_k\left\{\langle Y,Z_k\rangle + 2\frac{d_k}{r}\right\},
\nonumber
\end{eqnarray}
which implies
\begin{equation}\label{eqntilde4}
G_k(Y) = \frac{G(d_kY+X_k)}{d_k} \le (l+\e)\left\{\langle Y, Z_k\rangle +
2\frac{d_k}{r}\right\}.
\end{equation}
Passing to the limit as $k$ goes to infinity we conclude that for $Y\in
B(0,2)\cap\{\langle Y,e\rangle\ge 0\}$ $G_\infty(Y)\le (l+\e)\langle
Y,e\rangle$ for every $\e>0$, thus $G_\infty(Y)\le l\langle Y,e\rangle$.
Moreover $G_\infty(e)=l$. The maximum principle guarantees that
$v_\infty(Y)=l\max\{\langle Y,e\rangle;0\}$ for $Y\in B(0,1)$.

If $h_k(X)=h(d_kX+X_k)$, for $\zeta\in C^\infty_c(B(1,0))$, $\zeta\ge 0$
\begin{equation}\label{eqntilde5}
\int_{\p\{G_k>0\}}\zeta h_kd\H^n = \int_{\RR^{n+1}}\nabla
G_k\cdot\nabla\zeta\mathop{\longrightarrow}\limits_{k\to\infty} -
\int_{\RR^{n+1}}\nabla G_\infty\cdot\nabla\zeta = \int_{\{\langle
Y,e\rangle=0\}} l\zeta d\H^n
\end{equation}
thus
\begin{equation}\label{eqntilde6}
\lim_{k\to\infty} \int_{\p\{G_k>0\}}\zeta h_kd\H^n =l \int_{\{\langle
Y,e\rangle=0\}} \zeta dH^n.
\end{equation}
On the other hand the divergence theorem ensures that 
\begin{equation}\label{eqntilde7}
\int_{\p\{G_k>0\}}\zeta d\H^n\ge \int_{\p\{G_k>0\}}\zeta
e\cdot\nu_kd\H^n=\int_{\{G_k>0\}}\div(\zeta e).
\end{equation}
Since
\begin{equation}\label{eqntilde8}
\int_{\{G_k>0\}} \div(\zeta e)\mathop{\longrightarrow}\limits_{k\to\infty}
\int_{\{G_\infty>0\}}\div(\zeta e)=\int_{\p\{G_\infty>0\}}\zeta
d\H^n=\int_{\langle Y,e\rangle=0}\zeta d\H^n,
\end{equation}
we have that
\begin{equation}\label{eqntilde9}
\lim_{k\to\infty}\int_{\p\{G_k>0\}}\zeta d\H^n\ge \int_{\{\langle
Y,e\rangle=0\}}\zeta d\H^n.
\end{equation}
Since by (\ref{eqn16}), $h\ge e^{-\e}$ $\H^n-$ a.e. $Q\in\po$, using
(\ref{eqntilde6})
and (\ref{eqntilde9}) we have
\begin{eqnarray}\label{eqntilde10}
\lim_{k\to\infty} \int_{\p\{G_k>0\}} h_k\zeta d\H^n & \ge & \lim_{k\to\infty}
\int_{\p\{G_k>0\}}e^{-\e}\zeta d\H^n \\
l\int_{\{\langle Y,e\rangle=0\}}\zeta d\H^n & \ge & e^{-\e}\int_{\{\langle
Y,e\rangle=0\}} \zeta d\H^n, \nonumber
\end{eqnarray}
for any $\zeta\in C^\infty_c(B(1,0))$, $\zeta\ge 0$. Therefore
(\ref{eqntilde10}) yields
\begin{equation}\label{eqntilde11}
l\ge e^{-\e}.
\end{equation}
\qed
\end{proof}

\noindent Combining (\ref{eqn9}) and (\ref{eqntilde1}) we obtain that
\begin{equation}\label{eqn10}
|\nabla G_2(P_2)| \ge \limsup_{\mathop{X\to P_2}\limits_{X\in\O}}
\frac{G(X)}{d(X,B)}\ge e^{-\e}.
\end{equation}
Thus
\begin{equation}\label{eqn11}
\frac{1}{\H^n(\p B(0,R_2))} = \frac{1}{\s_nR^n_2}\ge e^{-\e}
\end{equation}
which implies 
\begin{equation}\label{eqn12}
R_2\le (e^\e\s^{-1}_n)^{\frac{1}{n}}.
\end{equation}
We have proved the following lemma.

\begin{lemma}\label{lem2}
Assume that $\O\subset\RR^{n+1}$ satisfies conditions (\ref{eqn1.1}) and
(\ref{eqn1.2}) in
Theorem \ref{thm1} then
\begin{equation}\label{eqn13}
B(0,R_1)\subset\O\subset B(0,R_2)
\end{equation}
with
\begin{equation}\label{eqn14}
e^{-\e}\le \s_nR^n_1\le\s_nR^n_2\le e^\e.
\end{equation}
\end{lemma}

%% file: stab-fine.tex
\section{Fine Geometric Properties}

In this section we prove Theorem \ref{thm1}. For this purpose we first
introduce a local notion of flatness that involves the geometry of the
boundary at a point $Q_0$, the behavior of $G$ near $Q_0$ and the
oscillation of $\log h$ near this point (see Definition 7.1 in \cite{AC}).
We assume that $G$ is continuously extended to be identically $0$ outside
$\O$. Note that $G$ is then subharmonic in $\RR^{n+1}$.

\begin{defn}\label{defn3.1}
Let $\O\subset\RR^{n+1}$ be as in Theorem \ref{thm1}. Let $Q_0\in\po$,
$\rho>0$ and $\s_+,\s_-, \tau\in (0,1)$. We say that
\begin{equation}\label{eqn3.1}
G\in F(\s_+,\s_-;\tau)\mbox{ in }B(Q_0,\rho)\mbox{ in direction
}\nu\mbox{ if}
\end{equation}
\begin{equation}\label{eqn3.2}
G(X)=0\mbox{ for }\langle X-Q_0, \nu\rangle\ge\s_+\rho 
\end{equation}
\begin{equation}\label{eqn3.3}
G(X) \ge -h(Q_0)[\langle X-Q_0,\nu\rangle+\s_-\rho]\mbox{ for }\langle
X-Q_0; \nu\rangle \le -\s_-\rho  
\end{equation}
and
\begin{equation}\label{eqn3.4}
\sup_{X\in B(Q_0,\rho)}|\nabla G(X)| \le h(Q_0)(1+\tau)\mbox{ and
}\mathop{\osc h}\limits_{B(Q_0,\rho)}\le\tau h(Q_0). 
\end{equation}
\end{defn}

The proof is very similar to the ones presented in \cite{AC} section 7 or in
\cite{KT1}. To avoid repetition we state the lemmata and only point out the
main differences with respect to the proofs of the results mentioned above.
For the complete details we refer the reader to \cite{AC} and \cite{KT1}.

\begin{lemma}\label{lem3.1}
Let $\O\subset\RR^{n+1}$ be a bounded domain and a set of locally finite
perimeter such that $0\in\O$. Let $G$ and $h$ be as above. There exists
$\s_n>0$ so that if $\s\in(0,\s_n)$, $\tau\in(0,\s)$ and $\e\in(0,\s)$ with 
\begin{equation}\label{eqn3.4a}
\sup_{\po}|\log h|<\e
\end{equation}
then for $Q_0\in\po$, $\rho>0$ and $\nu\in \SS^n$, if $G\in F(\s, 1;\tau)$ in
$B(Q_0,\rho)$ in direction $\nu$ then $G\in F(2\s, C\s;\tau)$ in 
$B\left(Q_0,\frac{\rho}{2}\right)$ in direction $\nu$. Here $C>1$ 
is a constant that only depends on $n$.
\end{lemma}

\begin{lemma}\label{lem3.2}
Let $\O\subset\RR^{n+1}$ be a bounded domain and a set of locally finite
perimeter such that $0\in \O$. Let $G$ and $h$ be as above. Given $\theta\in
(0,1)$ there exists $\s_\theta>0$ and $\eta_\theta=\eta\in(0,1)$ so that if
$\s\in(0,\s_\theta)$ and $\tau\in (0,\s_\theta\s^2)$ then for 
$Q_0\in\po$, $\rho>0$
if $G\in F(\s,\s;\tau)$ in $B(Q_0,\rho)$ in direction $\nu$ then $G\in
F(\theta\s, 1; \tau)$ in $B(Q_0,\eta\rho)$ in direction $\bar\nu$ and
$|\nu-\bar\nu|\le C\s$.
\end{lemma}

\begin{lemma}\label{lem3.3}
Assume that $\O\subset\RR^{n+1}$ satisfies (\ref{eqn1}). Then given $\s>0$
there exist $\e_\s>0$ such that if
\begin{equation}\label{eqn3.5}
\sup_{\po}|\log h|<\e\mbox{ with }\e<\e_\s
\end{equation}
then there is $\rho_\e=\rho>0$ (depending on $\e>0$) so that for $Q\in\po$,
$G\in F(\s,\s;(e^{2\e}-1)^{1/4})$ in $B(Q,\rho)$. Here
$(e^{2\e_\s}-1)^{1/2}<\s$.
\end{lemma}

\noindent{\bf Proof of Lemma \ref{lem3.3}}
Recall from Lemma \ref{lem2} that under the above hypothesis
$B(0,R_1)\subset\O\subset B(0, R_2)$ with $1\le R_2/R_1\le e^{2\e}$, and
$e^{-\e}\le\s_nR^n_i\le e^{\e}$ for $i=1,2$. Let
$\e\in\left(0,\frac{1}{4}\right)$ be a positive number to be chosen later
depending on $\s>0$. Let $\rho=R_1\sqrt{2}\sqrt{e^{2\e}-1}$. From basic
geometry and the remark above (see Lemma \ref{lem2}) it is clear that for
$Q\in\po$ there exists an $n$-plane $L(Q,\rho)$ through $Q$ such that
\begin{equation}\label{eqn3.6}
\frac{1}{\rho} D[\po\cap B(Q,\rho), L(Q,\rho)\cap B(Q,\rho)] \le \sqrt{2}
\sqrt{e^{2\e}-1}.
\end{equation}

In fact take for example the $n$-plane through $Q$ orthogonal to the line
joining the origin to $Q$. Let $\nu$ be the unit normal in the direction
$\overrightarrow{OQ}$ 
we have that if $X\in B(Q,\rho)$ and $\langle X-Q, \nu\rangle\ge
2\sqrt{2} \sqrt{e^{2\e}-1}\rho$ then since $1\le R_2/R_1\le e^{2\e}$
\begin{eqnarray}\label{eqn3.7}
|X|^2 & = & \left|X-Q-\langle X-Q,\nu\rangle\nu|^2+\langle X-Q,\nu\rangle+|Q|
\right|^2 \\
& \ge & \left(R_1 + 2\sqrt{2} \sqrt{e^{2\e}-1}\rho\right)^2 = R^2_1
\left(1+4(e^{2\e}-1)\right)^2 \nonumber \\
& \ge & R^2_2 e^{-4\e}(4e^{2\e}-3)^2 \nonumber \\
& \ge & R^2_2(4-3e^{-2\e})^2>R^2_2. \nonumber
\end{eqnarray}

Thus $X\not\in\O$ and $G(X)=0$ as $G$ was extended to be identically equal to
zero in $\O^c$. Now let $X\in B(Q,\rho)$ with $\langle X-Q,\nu\rangle \le
-2\sqrt{2} \sqrt{e^{2\e}-1}\rho$. In this case
\begin{eqnarray}\label{eqn3.8}
|X|^2 & = & \left|X-Q-\langle X-Q,\nu\rangle\nu|^2 + |\langle X-Q,\nu\rangle + 
|Q|\right|^2 \\
& \le & |X-Q|^2 + \left(R_2-4(e^{2\e}-1)R_1\right)^2 \nonumber \\
& \le & \rho^2+(e^{2\e}-4e^{2\e}+4)^2R^2_1 \nonumber \\
& \le & \left[2(e^{2\e}-1) + \left(1-3(e^{2\e}-1)\right)^2\right]R^2_1
\nonumber \\
& \le & \left(1+9(e^{2\e}-1)^2 - 4(e^{2\e}-1)\right)R^2_1 \nonumber \\
& \le & \left[1+(e^{2\e}-1)(9(e^{2\e}-1)-4)\right]R^2_1<R^2_1, \nonumber
\end{eqnarray}
provided $\e>0$ is such that $e^{2\e}-1<4/9$. Thus for $X\in B(Q,\rho)$ with
$\langle X-Q,\nu\rangle\le -2\sqrt{2}\sqrt{e^{2\e}-1}\rho$, $X\in B(0, R_1)$
and by (\ref{eqn5}) we have that if $G_1$ denotes the Green function of
$B(0, R_1)$ with pole $0$ then
\begin{eqnarray}\label{eqn3.9}
G(X) & \ge & G_1(X) = G_1(X)-G_1\left(R_1 \frac{X}{|X|}\right) \\
& \ge & -\sup_{Y\in B(Q,\rho)\cap B(0, R_1)} |\nabla G_1(Y)|(R_1-|X|).
\nonumber
\end{eqnarray}

The last inequality is a simple application of the fundamental theorem of
calculus. To estimate $\sup_{Y\in B(Q,\rho)\cap B(0,R_1)} |\nabla G_1(Y)|$
recall that $V_1(Y) = |\nabla G_1(Y)|$ is a subharmonic function on
$B(0,R_1)\backslash B\left(0, \frac{R_1}{2}\right)$. Since $R\le |Q|\le R_2$
then $B(Q,\rho)\subset B(0-, R_2+\rho)\backslash B(0, R_1-\rho)$. By the
maximum principle for bounded subharmonic functions 
\begin{equation}\label{eqn3.10}
\sup_{Y\in B(Q,\rho)\cap B(0,R_1)}|\nabla G_1(Y)| \le \sup_{Y\in\p
B(0,R_1)\cup\p B(0, R_1-\rho)}|\nabla G_1(Y)|.
\end{equation}
Since
$G_1(Y)=\frac{1}{(n-1)\s_n}
\left(\frac{1}{|Y|^{n-1}}-\frac{1}{R^{n-1}_1}\right)$, 
$\nabla G_1(Y)=\frac{-1}{\s_n}\frac{Y}{|Y|^{n+1}}$ and for $Y\in\p B(0,R_1)$
(\ref{eqn14}) implies that
\begin{equation}\label{eqn3.11}
|\nabla G_1(Y)| = \frac{1}{\s_nR^n_1} \le e^\e.
\end{equation}
For $Y\in\p B(0, R_1-\rho)$, our choice of $\rho$, and (\ref{eqn14}) ensure
\begin{eqnarray}\label{eqn3.12}
|\nabla G_1(Y)| & = & \frac{1}{\s_n} \frac{1}{(R_1-\rho)^n} \\
& = & \frac{1}{\s_nR^n_1} \frac{1}{(1-2(e^{2\e}-1))^n} \nonumber \\
& \le & \frac{e^\e}{(1-2(e^{2\e}-1))^n}. \nonumber
\end{eqnarray}

Combining (\ref{eqn3.9}), (\ref{eqn3.10}), (\ref{eqn3.11}) and
(\ref{eqn3.12}) we obtain
\begin{equation}\label{eqn3.13}
G(X)\ge e^\e (R_1-|X|)\mbox{ for }X\in B(Q,\rho)\mbox{ with }\langle X-Q,
\nu\rangle \le -2\sqrt{2} \sqrt{e^{2\e}-1}\rho.
\end{equation}

Our next goal is to compare $R_1-|X|$ to $|\langle X-Q, \nu\rangle|$. Note
that the basic picture is as follows:
\begin{center}
\setlength{\unitlength}{1cm}
\begin{picture}(8,8)
\put(4,3){\bigcircle{4}}
\put(4,3){\bigcircle{7.5}}
\put(4,6){\bigcircle[5]{3}}
\put(3.8,6){\line(0,-1){.2}}
\put(3.8,5.8){\line(1,0){.2}}
\put(4,3){\vector(0,1){5}}
\put(4,3){\line(-2,-3){1.1}}
\put(4,3){\line(-1,-5){.74}}
\put(-1,6){\line(1,0){10}}
\dottedline{.1}(4,3)(5.1,4.6)
\dottedline{.1}(4.39,4.9)(4.6,6)
\put(4.6,6){\circle*{.1}}
\put(3.26,-.655){\circle*{.1}}
\put(4,3){\line(1,5){.35}}
\put(4.62,6.4){\makebox(0,0)[t]{$\scriptstyle{\widehat X}$}}
\put(4.1,3.1){\makebox(0,0)[t]{$\scriptstyle{0}$}}
\put(3.3,1.6){\makebox(0,0)[t]{$\scriptstyle{R_1}$}}
\put(3.5,-.4){\makebox(0,0)[t]{$\scriptstyle{R_2}$}}
\put(4.9,5.6){\makebox(0,0)[t]{\begin{rotate}{90}\hbox{\Large{$\biggl\{$}}%
\end{rotate}}}
\put(4.7,5.6){\makebox(0,0)[t]{$\scriptstyle{\rho}$}}
\put(4.5,7){\makebox(0,0)[t]{$\scriptstyle{\nu}$}}
\put(4.2,6.4){\makebox(0,0)[t]{$\scriptstyle{Q}$}}
\put(4.5,4.8){\makebox(0,0)[t]{$\scriptstyle{X}$}}
\put(4.6,4.4){\makebox(0,0)[t]{$\scriptstyle{\theta}$}}
\put(4.45,4){\makebox(0,0)[t]{$\scriptstyle{\theta_X}$}}
\qbezier(4.25,4)(4.5,4.25)(4.75,4)
\qbezier(4,3.6)(4.25,3.85)(4.5,3.6)
\end{picture}
\end{center}
where
\begin{equation}\label{eqn3.14}
|\langle X-Q,\nu\rangle|=\cos\theta_X|X-\widehat X|\le
\left(R_1-|X|+\frac{R_2-R_1}{\cos\theta_X}\right)\cos\theta_X.
\end{equation}
Thus since $\langle X-Q, \nu\rangle \le -2\sqrt{2}\sqrt{e^{2_\e}-1}\rho\le
0$ and $R_1\le R_2\le e^{2\e}R_1$
\begin{eqnarray}\label{eqn3.15}
R_1-|X| & \ge & \frac{1}{\cos\theta_X} |\langle X-Q,\nu\rangle|-(R_2-R_1) \\
 & \ge & \frac{1}{\cos\theta_X} |\langle X-Q, \nu\rangle|-R_1(e^{2\e}-1)
\nonumber \\
 & \ge & \frac{1}{\cos\theta_X} |\langle X-Q, \nu\rangle| -
\frac{\sqrt{e^{2\e}-1}}{\sqrt{2}}\rho \nonumber \\
& \ge & |\langle X-Q,\nu\rangle|-\frac{e^{2\e}-1}{\sqrt{2}}\rho.
\nonumber
\end{eqnarray}

Combining (\ref{eqn3.5}), (\ref{eqn3.13}) and (\ref{eqn3.15}) we have that
for $X\in B(Q,\rho)$ with
$\langle X-Q;\nu\rangle\le -2\sqrt{2}\sqrt{e^{2\e}-1}\rho$
\begin{eqnarray}\label{eqn3.16}
G(X) & \ge & h(Q)(R_1-|X|) \\
& \ge & h(Q)\left[|\langle X-Q,\nu\rangle| -
\frac{e^{2\e}-1}{\sqrt{2}}\rho\right] \nonumber\\
& \ge & h(Q)\left[-\langle X-Q,\nu\rangle -
\frac{e^{2\e}-1}{\sqrt{2}}\rho\right]. \nonumber
\end{eqnarray}
Thus choosing $\e>0$ so that
$2\sqrt{2}\sqrt{e^{2\e}-1}<(e^{2\e}-1)^{\frac{1}{12}}<\sigma$ and we have
that for $X\in B(Q,\rho)$
\begin{equation}\label{eqn3.17}
G(X)=0\mbox{ for }\langle X-Q,\nu\rangle\ge \sigma\rho
\end{equation}
\begin{equation}\label{eqn3.18}
G(X) \ge - h(X)[\langle X-Q,\nu\rangle+\sigma\rho]\mbox{ for }\langle
X-Q,\nu\rangle\le -\sigma\rho.
\end{equation}
Hypothesis (\ref{eqn3.5}) implies that for $P, Q\in\po$
\begin{equation}\label{eqn3.19}
e^{-2\e}\le \frac{h(P)}{h(Q)}\le e^{2\e}.
\end{equation}

Thus
\begin{equation}\label{eqn3.20}
\osc_{B(Q,P)}h \le (e^{2\e}-1)h(Q).
\end{equation}

To estimate $\sup_{B(Q,P)\cap\O}|\nabla G|$ recall that the function $V(X)=
|\nabla G(X)|$ is subharmonic and bounded on $\O\bs B\left(0,
\frac{R_1}{2}\right)$. Hence since $B(Q,\rho)\cap\O\subset\O\bs B(0,
R_1-2\rho)$ the maximum principle for subharmonic functions ensures that
\begin{equation}\label{eqn3.21}
\sup_{B(Q,\rho)\cap\O}|\nabla G| \le \sup_{\O\bs B(Q,R_1-2\rho)}|\nabla G| =
\max\left\{\limsup_{X\to\po}|\nabla G(X)|, \sup_{\p B(0,R_1-2\rho)}|\nabla
G|\right\}.
\end{equation}

Let $Y\in\p B(0, R_1 -2\rho)$ then $B(Y,\rho) \subset B(0, R_1)\subset\O$
since $G$ and $G_1$ are harmonic on $B(Y,\rho)$ Poisson's representation
formula yields for $X\in B(Y,\rho)$
\begin{equation}\label{eqn3.22}
G(X)=\frac{\rho^2-|X-Y|^2}{(n+1)\o_{n+1}\rho} \int_{\p B(Y,\rho)}
\frac{G(\zeta)}{|X-\zeta|^{n+1}}d\zeta.
\end{equation}
Differentiating the expression in (\ref{eqn3.22}) and applying the obtained
formula to $X=Y$ we obtain
\begin{eqnarray}\label{eqn3.23}
\nabla G(Y) & = & -\frac{\rho}{\o_{n+1}} \int_{\p B(Y,\rho)}
\frac{G(\zeta)}{|Y-\zeta|^{n+3}}(Y-\zeta)d\zeta \\
& = & -\frac{1}{\o_{n+1}\rho^{n+2}} \int_{\p B(Y,\rho)}
G(\zeta)(Y-\zeta)d\zeta.\nonumber
\end{eqnarray}

Thus if $G_i$ denotes the Green function of $B(0, R_i)$ for $i=1,2$ 
with pole $0$, we have
\begin{equation}\label{eqn3.24}
|\nabla G(Y)-\nabla G_1(Y)| \le \frac{\rho}{\o_{n+1}\rho^{n+2}} \int_{\p
B(Y,\rho)} |G(\zeta)-G_1(\zeta)|d\zeta.
\end{equation}
Using (\ref{eqn5}), (\ref{eqn9}) and (\ref{eqn3.24}) we have
\begin{eqnarray}\label{eqn3.25}
|\nabla G(Y) - \nabla G_1(Y)| & \le & \frac{1}{\o_{n+1}\rho^{n+1}} \int_{\p
B(Y,\rho)} (G_2(\zeta)-G_1(\zeta))d\zeta \\
& \le & C_n \frac{1}{\rho^{n+1}} \int_{\p B(Y,\rho)} \left(
\frac{1}{R_1^{n-1}} - \frac{1}{R^{n-1}_2}\right)d\zeta \nonumber \\
& \le & \frac{C_n}{\rho}\left[\frac{1}{R^{n-1}_1} - \frac{1}{R^{n-1}_2}
\right] = \frac{C_n}{\rho R^{n-1}_1R^{n-1}_2} (R^{n-1}_2 - R^{n-1}_1)
\nonumber \\
& \le & \frac{C_nR^{n-2}_2}{\rho R^{n-1}_1R^{n-1}_2} (R_2-R_1) =
\frac{C_n}{\rho R^{n-1}_1R_2} (R_2-R_1) \nonumber \\
& \le & \frac{C_nR_1(e^{2\e}-1)}{\rho R^n_1} \le \frac{C_n}{\rho}
\frac{1}{R_1^{n-1}}(e^{2\e}-1) \nonumber \\
& \le & \frac{C_n}{R^n_1} \sqrt{e^{2\e}-1} \le C_n\sqrt{e^{2\e}-1},
\nonumber
\end{eqnarray}
where we used the facts that $1\le \frac{R_2}{R_1}\le e^{2\e}$,
$\rho=\sqrt{2}\sqrt{e^{2\e}-1}R_1$ and $e^{-\e}\le R^n_1\s_n\le e^\e$, with
$\e\in\left(0,\frac{1}{4}\right)$. Since
$G_1(Y)=\frac{1}{(n-1)(n+1)\o_{n+1}} \left(\frac{1}{|Y|^{n-1}} -
\frac{1}{R^{n-1}_1}\right)$ then $|\nabla G_1(Y)| = \frac{1}{(n+1)\o_{n+1}}
\frac{1}{|Y|^n}$. For $Y\in\p B(0,R_1-2\rho)$ and $\e>0$ small enough, we
have 
\begin{eqnarray}\label{eqn3.26}
|\nabla G_1(Y)| & = & \frac{1}{\o_{n+1}(n+1)(R_1-2\rho)^n} \\
& = & \frac{1}{\s_nR^n_1(1-2\sqrt{2}\sqrt{e^{2\e}-1)}^n} \nonumber \\
& \le & \frac{e^\e}{(1-2\sqrt{2}\sqrt{e^{2\e}-1})^n} \nonumber \\
& \le & e^\e(1+8n\sqrt{e^{2\e}-1}). \nonumber
\end{eqnarray}

Combining (\ref{eqn3.5}), (\ref{eqn3.21}), (\ref{eqn3.25}) and
(\ref{eqn3.26}) we obtain 
\begin{eqnarray}\label{eqn3.27}
\sup_{B(Q,\rho)}|\nabla G| & \le & e^\e(1+8n\sqrt{e^{2\e}-1}) +
C_n\sqrt{e^{2\e}-1} \\
& \le & e^\e(1+C_n\sqrt{e^{2\e}-1}) \nonumber \\
& \le & e^{2\e}h(Q)(1+C_n\sqrt{e^{2\e}-1}) \nonumber \\
& \le & h(Q)(1+C_n\sqrt{e^{2\e}-1}). \nonumber
\end{eqnarray}

Thus for $\e>0$ small enough so that $C_n(e^{2\e}-1)^{\frac{1}{4}}<1$ we have
that
\begin{equation}\label{eqn3.28}
\sup_{B(Q,\rho)} |\nabla G| \le h(Q)(1+(e^{2\e}-1)^{1/4}).
\end{equation}
Note that (\ref{eqn3.17}), (\ref{eqn3.18}), (\ref{eqn3.20}) and
(\ref{eqn3.28}) show that for $\e>0$ small enough in terms of $n$ and such
that $(\e^{2\e}-1)^{1/12}<\s$ then $G\in F(\s,\s;(e^{2\e}-1)^{1/4})$ in
$B(Q,\rho)$, $\forall\,Q\in\po$ where $\rho=\sqrt{2}\sqrt{e^{2\e}-1}R_1$.
\qed

Before sketching the proofs of Lemma \ref{lem3.1} and Lemma \ref{lem3.2} we
indicate how from the 3 lemmata above one proves Theorem \ref{thm1}.
\medskip

\noindent{\bf Proof of Theorem \ref{thm1}:} Let $\theta'\in\left(0,
\frac{1}{2}\right)$ to be chosen. Let $\s'\in(0, \s_{\theta'})$ as in Lemma
\ref{lem3.2}. By Lemma \ref{lem3.3} for $\s\in (0, \s_{\theta'})$ there
is $\e_{\s'}>0$ so that if (\ref{eqn3.5}) holds, then 
$G\in F(\s', \s', (e^{2\e}-1)^{1/4})$ in
$B(Q,\rho)$, for $Q\in\po$ with
$\rho=\sqrt{2}\sqrt{e^{2\e}-1}R_1$, and with $(e^{2\e_{\s'}}-1)^{1/12}<\s'$.
Note that by choosing $\e'<\e_{\s'}$ so that
$(e^{2\e_\s}-1)^{1/4}<\s_{\theta'}$ we have that
$(e^{2\e}-1)^{1/4}\le \s_{\theta'}(\s')^2$ for $\e<\e'$. Lemma \ref{lem3.2}
ensures that $G\in F(\theta'\s', 1; (e^{2\e}-1)^{1/4})$ in $B(Q,\eta\rho)$.
Lemma \ref{lem3.1} now guarantees that $G\in F(2\theta'\s', C\theta'\s';
(e^{2\e}-1)^{1/4})$ in $B\left(Q, \frac{\eta\rho}{2}\right)$. Choosing
$\theta'$ so that $C\theta'+2\theta'<1$ we conclude that $G\in F(\s', \s';
(e^{2\e}-1)^{1/4})$ in $B\left(Q, \frac{\eta\rho}{2}\right)$. Since
$(e^{2\e}-1)^{1/4}\le\s_{\theta'}(\s')^2$ we can repeat the previous
argument to show that $\forall\,k\in\NN$ and $\forall\,Q\in\po$ $G\in
F(\s',\s'; (e^{2\e}-1)^{1/4})$ in $B\left(Q,
\left(\frac{\eta}{2}\right)^k\rho\right)$.

Thus there exists $\nu_k\in S^n$ so that
\begin{equation}\label{eqn3.29}
G(X)=0\mbox{ for }\langle X-Q,\nu_k\rangle\ge
\s'\left(\frac{\eta}{2}\right)^k\rho
\end{equation}
and
\begin{equation}\label{eqn3.30}
G(X) \ge -h(Q)\left[\langle X-Q, \nu_k\rangle +
\s'\left(\frac{\eta}{2}\right)^k\rho\right]\ge 0 \mbox{ for }\langle
X-Q,\nu_k\rangle \le -\s'\left(\frac{\eta}{2}\right)^k\rho.
\end{equation}

In particular if $L_k(Q)$ denotes the $n$-plane through $Q$ orthogonal to
$\nu_k$ (\ref{eqn3.29}) and (\ref{eqn3.30}) imply that
\begin{equation}\label{eqn3.31}
D\left[\po\cap B\left(Q, \left(\frac{\eta}{2}\right)^k\rho\right);
L_k(Q)\cap B\left(Q, \left(\frac{\eta}{2}\right)^k\rho\right)\right] \le
\s'\left(\frac{\eta}{2}\right)^k\rho.
\end{equation}
Let $r\in (0, \rho)$ there is $k\ge 0$ so that
$\left(\frac{\eta}{2}\right)^{k+1}\rho\le r\le
\left(\frac{\eta}{2}\right)^k\rho$, let
$r_k=\left(\frac{\eta}{2}\right)^k\rho$. For $P\in\po\cap B(Q,r)$ by
(\ref{eqn3.31}), there exists $Z\in L_k(Q)\cap B(Q,r_k)$ so that
$|Z-P|<\s'r_k$. Note that $|Z-Q|\le |Z-P|+|P-Q|<\s'r_k+r$. There exists
$Z'\in \seg[Q,Z]$ such that $|Z'-Q|<r$ and $|Z'-Z|<\s'r_k$. Moreover
$|Z'-P|\le |Z-Z'|+|Z-P|<2\s'r_k$. 

For $Z\in L_k(Q)\cap B(Q,r)$, there exists $Z'\in L_k(Q)\cap B(Q,r-\s'r_k)$
so that $|Z-Z'|<\s'r_k$. By (\ref{eqn3.31}) there exists $P\in\po\cap
B(Q,r_k)$ so that $|Z'-P|<\s'r_k$. Note that $|Z-P|\le
|Z-Z'|+|Z'-P|<2\s'r_k$, moreover $|P-Q|\le |P-Z'|+|Z'-Q|<r$. Thus
$P\in\po\cap B(Q,r)$. The previous argument ensures that for $Q\in\po$ and
$r\in (0,\rho)$ there exists an $n$-plane through $Q$, $L(Q,r)$ so that
\begin{equation}\label{eqn3.32}
\frac{1}{r}D[\po\cap B(Q,r), L(Q,r)\cap B(Q,r)]\le 2\s'.
\end{equation}

Thus for $\s\in\left(0, \frac{\s_{\theta'}}{2}\right)$ there exists
$\e_\s>0$ so that if $\e<\e_\s$ and $\sup_{\po}|\log h|<\e$ then
$\theta(Q,r)\le\s$ for $r\in (0; \rho)$ with
$\rho=\sqrt{2}\sqrt{e^{2\e}-1}R_1$.
\qed
\medskip

We now focus our attention in the proofs of Lemmas \ref{lem3.1} and
\ref{lem3.2}. As mentioned earlier these are just small variations of
results that appear both in \cite{AC} and \cite{KT1}, thus we do not present
all the details.
\medskip

\noindent{\bf Proof of Lemma \ref{lem3.1}:} Without loss of generality we
may assume that $Q_0=0\in\po$, $\rho=1$ and $\nu=e_{n+1}$. By hypothesis $G\in
F(\s,1;\tau)$ in $B_1=B(0,1)$ in the direction $e_{n+1}$, $h(Q)\ge e^{-\e}$
for $\H^n$ a.e. $Q\in\po$ and $\sup_{B_1}|\nabla G|\le e^\e(1+\tau)\le
e^\e(1+\s)$. This implies that for $\varphi\in C^\infty_0(\RR^{n+1})$,
$\varphi\ge 0$
\begin{equation}\label{eqn3.34}
\int_\O G\Delta\varphi\ge e^{-\e}\int_{\po}\varphi d\H^n.
\end{equation}
Let $\eta(Y)=\exp\left(\frac{-9|Y|^2}{1-9|Y|^2}\right)$ for
$|Y|<\frac{1}{3}$ and $\eta(Y)=0$ otherwise. Choose $s_0>0$ to be the
maximum $s$ so that
\begin{equation}\label{eqn3.35}
B_1\cap\{G>0\}\subset D=\{X\in B_1:x_{n+1}<2\s-s\eta(\bar x)\}
\end{equation}
where $X=(\bar x, x_{n+1})$ with $\bar x\in\RR^n\times\{0\}$ Note that
$s_0\le 2\s$. Since $G\in F(\s,1;\tau)$ in $B_1$ there exists $Z\in\p D\cap
\po\cap B\left(0, \frac{1}{3}\right)$. Let $B\subset D^C$ be a tangent ball
to $D$ at $Z$. Since $\p D\cap B_1$ is smooth and $s_0\le 2\s\le\s_n$ for
$\s_n>0$ small we may assume that the radius of $B$ is $\frac{C_n}{\s_n}$.
Consider the function $V$ defined by $\Delta V=0$ in $D$, $V=0$ in $\p D\cap
B_1$ and $V=2\s-x_{n+1}$ on $\p D\bs B_1$. By the maximum principle $V>0$ in
$D$ and
\begin{equation}\label{eqn3.36}
G\le V\mbox{ in }D
\end{equation}
as $G\le V$ on $\p D$ and $G$ is subharmonic. For $X\in D$ define
$F(X)=(2\s-x_{n+1})-V(X)$, $F$ is a harmonic function on $D$. Since $Z$ is a
smooth point of $\p D$, standard boundary regularity arguments (see
\cite[Lemma 6.5]{GT}) ensures that $\sup_{X\in\bar D}|\nabla F(X)|\le
C\sup_{\bar D}|F|\le Cs_0\le C\s$. Therefore
\begin{equation}\label{eqn3.38}
-\frac{\p V}{\p x_{n+1}} (Z)=1 + \frac{\p F}{\p x_{n+1}}(Z)\le 1+C\s.
\end{equation}

Using (\ref{eqn3.38}) and noting that $|\vecn(Z)-e_{n+1}|\le c\s$ we have
that if $\langle\nabla V, \vecn\rangle=\frac{\p V}{\p n}$ where $\vecn$
denotes the outward unit normal to $\p D$ then
\begin{equation}\label{eqn3.39}
-\frac{\p V}{\p b}(Z)\le 1+c\s+(1+\s)|\vecn-e_{n+1}|\le 1+c\s.
\end{equation}
Our goal now is to estimate $G$ from below by the linear function $-x_{n+1}$
up to a constant of order $\s$. Let $\zeta\in\p B\left(0,
\frac{3}{4}\right)\cap\left\{x_{n+1}<-\frac{1}{2}\right\}$. Consider the
function $\o_\zeta$ defined by $\Delta\o_\zeta=0$ in $D\bs
B\left(\zeta,\frac{1}{8}\right)$, $\o_\zeta=0$ on $\p D$ $\o_\zeta=-x_{n+1}$
on $\p B\left(\zeta,\frac{1}{8}\right)$. The Hopf boundary point lemma
ensures that
\begin{equation}\label{eqn3.40}
-\frac{\p \o_\zeta}{\p n}(Z)\ge C_n>0.
\end{equation}
Assume that there exists $d>0$ such that $\forall\,X\in\bar B\left(\zeta,
\frac{1}{8}\right)$
\begin{equation}\label{eqn3.41}
G(X)\le V(X)+\s dx_{n+1}.
\end{equation}
The maximum principle would then imply that
\begin{equation}\label{eqn3.42}
G(X)\le V(X)-d\s\o_\zeta(X)\mbox{ in }D\bs B\left(\zeta,\frac{1}{8}\right).
\end{equation}

Combining Lemma \ref{lem3.1}, (\ref{eqn3.39}), (\ref{eqn3.32}),
(\ref{eqn3.4a}) and the hypothesis that $\e\in(0,\s)$ we would have
\begin{equation}\label{eqn3.43}
1-\s\le 1-\e  \le -\frac{\p V}{\p n}(Z)-d\s\frac{\p\o_\zeta}{\p n}(Z) 
 \le  1+C\s -C_nd\s 
\end{equation}
which is a contradiction for $d$ large. Thus fr $d$ large enough
(depending on $n$) there are points $X_\zeta\in
B\left(\zeta, \frac{1}{8}\right)$ such that
\begin{equation}\label{eqn3.44}
G(X_\zeta)\ge V(X_\zeta)+d\s(X_\zeta)_{n+1}.
\end{equation}

Let $X\in B\left(X_\zeta, \frac{1}{4}\right)$ then noting that $V(X)\ge
-x_{n+1}$ for $X\in D$, using the fact that $\sup_{B_1}|\nabla G|\le
e^\e(1+\s)$ and (\ref{eqn3.44}) we have for $\s_n$ small enough
\begin{eqnarray}\label{eqn3.45}
G(X) & \ge & G(X_\zeta) - \sup_{B\left(\zeta, \frac{1}{4}\right)} |\nabla
G|\,|X-X_\zeta| \\
& \ge & V(X_\zeta) + d\s(X_\zeta)_{n+1} - \frac{1}{4}(1+\s)e^\e \nonumber \\
& \ge & -(X_\zeta)_{n+1} + d\s(X_\zeta)_{n+1} - \frac{1}{4}(1+\s)e^\e
\nonumber \\
& \ge & \frac{5}{8} - \frac{7}{8}d\s - \frac{1}{4}(1+\s)e^\e \nonumber \\
& \ge & \frac{5}{8} - \frac{7}{8}d\s - \frac{1}{4}(1+\s)e^\s>0 \nonumber
\end{eqnarray}
for $\s<\s_n$. Since $G(X)>0$ for $X\in\overline{B\left(X_\zeta,
\frac{1}{4}\right)}$, $G$ is harmonic on $B\left(X_\zeta, \frac{1}{4}\right)$
and so is $V-G$. Moreover $V-G\ge 0$ on $B\left(X_\zeta,
\frac{1}{4}\right)\supset B\left(\zeta, \frac{1}{8}\right)$.

Harnack's inequality combined with (\ref{eqn3.44}) yields
\begin{equation}\label{eqn3.46}
(V-G)(\xi)  \le  C_n(V-G)(X_\zeta)  \le  -Cd\s(X_\zeta)_{n+1}\le C\s 
\end{equation}
and
\begin{equation}\label{eqn3.47}
G(\zeta)\ge V(\zeta)-C\s\ge -\zeta_{n+1}-C\s.
\end{equation}
For $X\in D\cap B\left(0, \frac{1}{2}\right)$, $X=\zeta+tx_{n+1}$ for some
$\zeta\in\p B\left(0, \frac{3}{4}\right)\cap
\left\{x_{n+1}<-\frac{1}{2}\right\}$ then (\ref{eqn3.47}) implies that
\begin{equation}\label{eqn3.48}
G(X)  \ge  G(\zeta)-(1+\s)e^\s t  \ge  -(\zeta_{n+1}+t)-C\s 
\end{equation}
since $G\in F(\s, 1;\tau)$ in $B_1$ in  direction $e_{n+1}$, inequality
(\ref{eqn3.48}) ensures that $G\in F(2\s, C\s; \tau)$ in $B\left(0,
\frac{1}{2}\right)$ in  direction $e_{n+1}$.\qed
\medskip

Lemma \ref{lem3.2} is proved by contradiction, using a non-homogeneous
blow-up. Assume that Lemma \ref{lem3.2} does not hold. There exists
$\theta_0\in (0,1)$ such that for every $\eta>0$ (later we specify one) and
every non-negative decreasing sequence $\{\s_j\}$ there is a sequence
$\{\tau_j\}$ with $\tau_j\s_j^{-2} \to 0$ so that
\begin{equation}\label{eqn3.49}
G\in F(\s_j, \s_j; \tau_j)\mbox{ in }B(Q_j, \rho_j)\mbox{ in direction
}\nu_j
\end{equation}
but
\begin{equation}\label{eqn3.50}
G\not\in F(\theta_0\s_j, 1; \tau_j)\mbox{ in }B(Q_j,\eta\rho_j).
\end{equation}

Since the estimate in Lemma \ref{lem3.2} is to hold uniformly on compact
sets we assume that for each $j\in \NN$, $Q_j\in K$ and that
$\lim_{j\to\infty}Q_j=Q_0\in K$ $Q_0\ne 0$ where $K$ is a fixed compact set
in $\RR^{n+1}$.

Note that if $G\in F(\s,\s;\tau)$ in $B(Q,\rho)$ in direction $\nu$ then
$G\in F(4\s, 4\s; \tau)$ in $B\left(P, \frac{\rho}{2}\right)$ in direction
$\nu$ for every $P\in\po\cap B\left(Q, \frac{P}{2}\right)$. Let $R_j$ be the
rotation which maps $\RR^{n+1}_+$ onto $\{(x,t)=x+t\nu_j:
x\in\langle\nu_j\rangle^{\perp}; t\ge 0\}$. Let $\O_j=\rho^{-1}_j
R^{-1}_j(\O-Q_j)$, $\po_j=\rho^{-1}_jR^{-1}_j(\po-Q_j)$. Define
\begin{equation}\label{eqn3.51}
G_j(X) = \frac{1}{\rho_jh(Q_j)} G(\rho_jR_jX+Q_j)
\end{equation}
and for $Q\in\po_j$
\begin{equation}\label{eqn3.52}
h_j(Q)=\frac{1}{h(Q_j)} h(\rho_jR_jQ+Q_j).
\end{equation}

Note that $G_j$ is a positive multiple of the Green function of $\O_j$ with
pole $-\rho_j^{-1} R^{-1}_jQ_j$. Note that $|\rho^{-1}_jR^{-1}_jQ_j|\ge
\frac{|Q_0|}{2\rho_j}$ for $j$ large enough. Thus for $\varphi\in
C^\infty_c(\RR^{n+1})$ and $j$ large enough so  ${\rm support}\,\varphi\subset
B\left(0, \frac{|Q_0|}{4\rho_j}\right)$ we have
\begin{equation}\label{eqn3.53}
\int_{\O_j}G_j\Delta\varphi dX = \int_{\po_j}\varphi h_jd\H^n
\end{equation}
with
\begin{equation}\label{eqn3.54}
\sup_{B(0,1)}|\nabla G_j| \le 1+\tau_j\mbox{ and
}\osc_{B(0,1)}h_j\le\tau_j\mbox{ with }h_j(0)=1.
\end{equation}
Moreover
\begin{equation}\label{eqn3.55}
G_j\in F(\s_j, \s_j; \tau_j)\mbox{ in }B(0,1)\hbox{ in direction }e_{n+1}
\end{equation}
but
\begin{equation}\label{eqn3.56}
G_j\not\in F(\theta_0\s_j, 1; \tau_j)\mbox{ in }B(0,\eta)
\end{equation}
with $\s_j\to 0$ and $\tau_j\s_j^{-2}\to 0$ as $j\to \infty$.

We define sequences of scaled height functions (in the direction $e_{n+1}$)
corresponding to $\po_j$. We prove that this sequence converges to a
subharmonic Lipschitz function, and use this information to contradict
(\ref{eqn3.56}) for $j$ large enough. For $y\in
B(0,1)\cap\RR^n\times\{0\}=B'$ define
\begin{equation}\label{eqn3.57}
f^+_j(Y)=\sup\left\{h:(y_1\s_jh)\in\p\{G_j>0\}\right\}\le 1
\end{equation}
and
\begin{equation}\label{eqn3.58}
f^-_j(Y) = \inf\left\{h:(y, \s_jh)\in\p\{G_j>0\}\right\}\ge -1
\end{equation}

\begin{lemma}[Non-homogeneous blow up (Lemma 7.3 \cite{AC})]\label{lem3.4}
There exists a subsequence $k_j$ such that for $y\in B'$
\begin{equation}\label{eqn3.59}
f(y) = \limsup_{\mathop{k_j\to\infty}\limits_{z\to y}} f^+_{k_j}(z) =
\liminf_{\mathop{k_j\to\infty}\limits_{z\to y}} f^-_{k_j}(z).
\end{equation}
\end{lemma}

\begin{corollary}[Corollary 7.4 \cite{AC}]\label{cor3.1}
The function $f$ that appears in (\ref{eqn3.59}) is a continuous function in
$B'$, $f(0)=0$; and $f^+_{k_j}$ and $f^-_{k_j}$ converge uniformly to $f$ on
compact sets of $B'$.
\end{corollary}

The proofs of Lemma \ref{lem3.4} and Corollary \ref{cor3.1} are identical to
those that appear in \cite{AC} or \cite{KT1}, thus we omit them here.

\begin{lemma}[Lemma 7.5 \cite{AC}]\label{lem3.5}
The function $f$ introduced in Lemma \ref{lem3.4} is subharmonic in $B'$.
\end{lemma}

\begin{proof}
This proof is done by contradiction. Assuming that $f$ is not subharmonic
in $B'$ we contradict the fact that $\s^{-2}_j\tau_j\to 0$ as $j\to\infty$.
In fact if $f$ is not subharmonic in $B'$ there exists $y_0\in B'$ and
$\rho>0$ so that $B'(y_0,\rho)\subset B'$ and
\begin{equation}\label{eqn3.60}
f(y_0) > \notint_{\lower.1in\hbox{$\scriptstyle{\p B'(y_0,\rho)}$}}f(x)dx.
\end{equation}
\end{proof}

Let 
\begin{equation}\label{eqn3.61}
\e_0=\frac{f(y_0)-\notint_{\p B'(y_0,\rho)}f(x)dx}{2}.
\end{equation}
Let $g$ be the solution to the Dirichlet problem
\begin{equation}\label{eqn3.62}
\left\{\begin{array}{rclll}
\Delta g & = & 0 & \mbox{in} & B'(y_0,\rho) \\
g & = & f+\e_0 & \mbox{on} & \p B'(y_0,\rho). 
\end{array}
\right\}
\end{equation}
Note that
\begin{equation}\label{eqn3.63}
f<g\mbox{ on }\p B'(y_0,\rho),\mbox{ and } 
\end{equation}
\begin{eqnarray}\label{eqn3.64}
g(y_0) & = & \notint_{\p B'(y_0,\rho)}g(x)dx = \notint_{\p
B'(y_0,\rho)}f(x)dx+\e_0 \\
g(y_0) & = & \frac{1}{2}\left\{f(y_0) + \notint_{\p B'(y_0,\rho)}f(x)dx
\right\} \nonumber \\
g(y_0) & < & f(y_0).
\end{eqnarray}
Summarizing, we have the following picture.

\begin{center}
\begin{psfrags}
\includegraphics[width=3in]{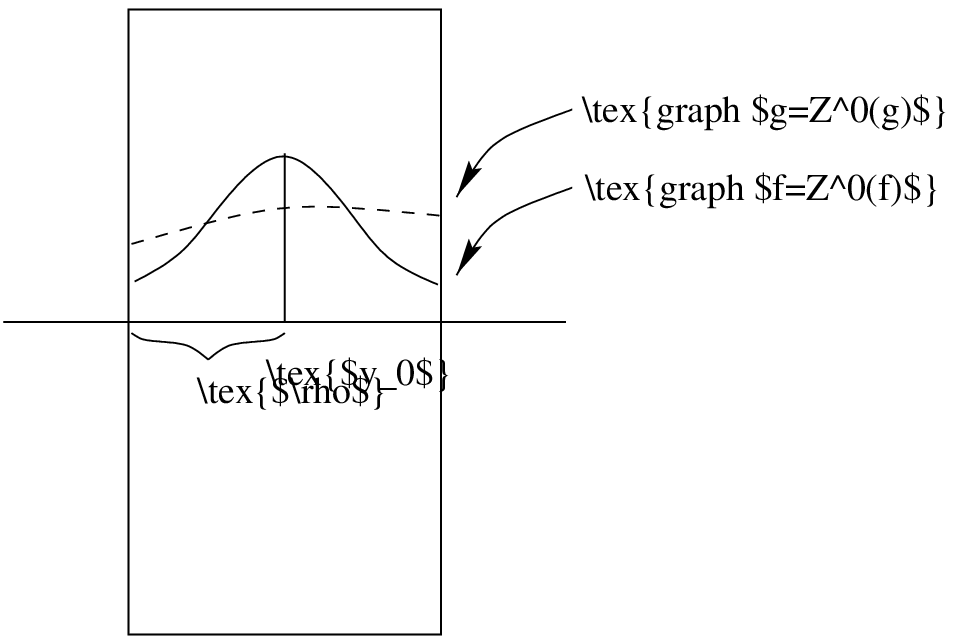}
\end{psfrags}
\end{center}

\begin{equation}\label{eqn3.65}
\left\{\begin{array}{rclll}
\Delta g & = & 0 & \mbox{in} & B'(y_0,\rho) \\
g & > & f & \mbox{in} & \p B'(y_0,\rho) \\
g(y_0) & < & f(y_0) & 
\end{array}\right.
\end{equation}
The main idea of the proof is to compare the $n$-dimensional Hausdorff
measure of $\p\{G_{k_j}>0\}$ on the cylinder $B'(y_0,\rho)\times(-1,1)$ to
that of the graph of $\s_{k_j}g$ on the same cylinder to obtain a
contradiction from an estimate on the size of the area enclosed by these 2
surfaces. In order to simplify the notation we relabel the sequences that
appear in Lemma \ref{lem3.4}. We also introduce some new definitions.

Let $Z=B'(y_0,\rho)\times\RR$ be the infinite cylinder. For $\phi$ defined
on $\RR^n$ define
\begin{eqnarray}\label{eqn3.66}
Z^+(\phi) & = & \{(y,h)\in Z:h>\phi(y)\} \\
Z^-(\phi) & = & \{(y,h)\in Z: h<\phi(y)\} \nonumber \\
Z^0(\phi) & = & \{(y,h)\in Z:h=\phi(y)\}. \nonumber
\end{eqnarray}
We may assume that for $k$ large enough
\begin{equation}\label{eqn3.67}
\H^n(Z^0(\s_kg)\cap\p\{G_k>0\})=0.
\end{equation}
(It might be necessary to modify $g$ above by adding a suitable constant
which can be chosen as small as one wants. In particular the function $g$
would still satisfy (\ref{eqn3.63}) and (\ref{eqn3.64}).

\begin{Claim}\label{claim1}
For $k$ large enough
\begin{equation}\label{eqn3.68}
\H^n(Z^+(\s_kg)\cap\p\{G_k>0\}) \le \frac{1+\tau_k}{1-\tau_k}
\H^n(Z^0(\s_kg)\cap \{G_k>0\}).
\end{equation}
\end{Claim}

\begin{Claim}\label{claim2}
Let $E_k=\{G_k>0\}\cap Z^-(\s_kg)$. $E_k$ is a set of locally finite
perimeter and
\begin{equation}\label{eqn3.69}
\H^n(Z\cap\p^\ast E_k) \le \H^n(\p\{G_k>0\}\cap Z^+(\s_kg)) + \H^n(\{G_k=0\}\cap
Z^0(\s_kg)).
\end{equation}
Here $\p^\ast E_k$ denotes the reduced boundary of $E_k$.
\end{Claim}

\begin{Claim}\label{claim3}
There exists a constant $C>0$ such that 
\begin{equation}\label{eqn3.70}
\H^n(Z\cap\p^\ast E_k)\ge\H^n(Z^0(\s_kg))+C\s^2_k\rho^n.
\end{equation}
\end{Claim}

Before proving the claims we indicate how combining inequalities
(\ref{eqn3.68}), (\ref{eqn3.69}) and (\ref{eqn3.70}) we obtain a
contradiction. Combining (\ref{eqn3.68}), (\ref{eqn3.69}) and (\ref{eqn3.70}) and
using (\ref{eqn3.67}) we have
\begin{eqnarray}\label{eqn3.71}
&&\kern-.5in\H^n(Z^0(\s_kg))+C\s^2_k\rho^n  \le  \H^n(Z\cap\p^\ast E_k) \\
& \le & \H^n(\p\{G_k>0\}\cap Z^+(\s_kg)+\H^n(\{G_k=0\}\cap Z^-(\s_kg))
\nonumber \\
& \le & \frac{1+\tau_k}{1-\tau_k}\H^n(Z^0(\s_kg)\cap\{G_k>0\}) +
\H^n(\{G_k=0\}\cap Z^0(\s_kg)) \nonumber \\
& \le & \frac{2\tau_k}{1-\tau_k}\H^n(Z^0(\s_kg)\cap\{G_k>0\}) + \H^n(Z^0(\s_kg))
\nonumber
\end{eqnarray}
which implies
\begin{eqnarray}\label{eqn3.72}
C\s^2_k\rho^n  & \le & \frac{2\tau_k}{1-\tau_k}\H^n(Z^0(\s_kg)\cap\{G_k>0\})
\\
& \le & \frac{2\tau_k}{1-\tau_k}\int_{B_{\rho'(y_0)}}\sqrt{1+\s^2_k|\nabla
g|^2}. \nonumber
\end{eqnarray}

For $\tau_k<\frac{1}{2}$ and $\s_k<1$ (\ref{eqn3.72}) yields $C\s^2_k\le
C'\tau_k$ which contradicts the fact that $\tau_k\s^{-2}_k\to 0$ as
$k\to\infty$. Thus we conclude that $f$ is subharmonic in $B'$.
\medskip

\noindent{\bf Proof of Claim \ref{claim1}:} Since $h_k(0)=1$ and
$\osc_{B(0,1)}h_k\le\tau_k$ we have that
\begin{eqnarray}\label{eqn3.73}
\H^n(Z^+(\s_kg)\cap\p\{G_k>0\}) & = & \int_{Z^+(\s_kg)\cap\p\{G_k>0\}}d\H^n
\\
& \le & \frac{1}{1-\tau_k}
\int_{Z^+(\s_kg)\cap\p\{G_k>0\}}h_kd\H^n.\nonumber
\end{eqnarray}
For $\varphi\in C^\infty_c(\RR^{n+1})$ and $k$ large enough we have
\begin{equation}\label{eqn3.74}
-\int_{\{G_k>0\}}\nabla G_k\nabla\varphi = \int_{\p\{G_k>0\}}\varphi
h_kd\H^n.
\end{equation}
Letting $\varphi\to\chi_{Z^+(\s_kg)}$, (\ref{eqn3.74}) yields
\begin{equation}\label{eqn3.75}
-\int_{\{\overline{G_k>0}\}\cap\p Z^+(\s_kg)}\nabla G_k\cdot\nu =
\int_{\p\{G_k>0\}\cap Z^+(\s_kg)}h_kd\H^n
\end{equation}
where $\nu$ denotes the outward pointing unit normal. Combining
(\ref{eqn3.67}), (\ref{eqn3.73}), (\ref{eqn3.75}) and (\ref{eqn3.54}) we
have that
\begin{eqnarray}\label{eqn3.76}
\H^n(Z^+(\s_kg)\cap\p\{G_k>0\}) & \le & \frac{1}{1-\tau_k}
\int_{\{G_k>0\}\cap\p Z^+(\s_kg)} |\nabla G_k| \\
& \le & \frac{1+\tau_k}{1-\tau_k} \H^n(\{G_k>0\}\cap Z^0(\s_kg)).\nonumber
\end{eqnarray}
\qed
\medskip

The proof of Claim \ref{claim2} is straightforward. The proof of Claim
\ref{claim3} is identical to the one that appears in either \cite{AC} or
\cite{KT1}, thus we do not present it here.

To obtain the desired contradiction we need to prove that $f$ is Lipschitz. 
This proof relies on the following lemma
which claims that $f$ converges to its average faster than linearly in an
integral sense.

\begin{lemma}[Lemma 7.6 \cite{AC}]\label{lem3.6}
There is a constant $C=C(n)>0$ such that for $y\in B'_{1/2}=B\left(0,
\frac{1}{2}\right)\cap\RR^n\times\{0\}$
\begin{equation}\label{eqn3.61b}
0\le\int^{\frac{1}{4}}_0 \frac{1}{r^2} (f_{y,r} - f(y))dr\le C
\end{equation}
where
\begin{equation}\label{eqn3.62b}
f_{y,r} = \notint_{\lower.1in\hbox{$\scriptstyle{\p B'(y,r)}$}}fd\H^{n-1}.
\end{equation}
\end{lemma}

\begin{proof}
The proof is very similar to the ones that appear in \cite{AC} and
\cite{KT1}. Nevertheless since the minor differences are technically
important we sketch the proof here pointing out how to overcome the
difficulties that arise in this situation. For the complete details we refer
the reader to \cite{AC} or \cite{KT1}. Without loss of generality we may
assume that $y=0$. Since $f(0)=0$ it is enough to show
\begin{equation}\label{eqn3.63b}
0\le \int^{\frac{1}{4}}_0 \frac{1}{r^2} 
\notint_{\lower.1in\hbox{$\scriptstyle{\p B'_r}$}} fd\H^{n-1}\le C
\end{equation}
where $B'_r=B'(0,r)$ and $C$ only depends on $n$ since $f$ is subharmonic
(see Lemma \ref{lem3.5}) then for $r\in \left(0, \frac{1}{2}\right)$, 
$f(0)\le\notint_{\p B'_r}fd\H^{n-1}$ which proves the first inequality.

Let $h>2\s_j$ be small and let $G_h$ denote the Green function of $B\left(0,
\frac{1}{2}\right)\cap\{x_{n+1}<0\}$ with pole $-he_{n+1}$. By reflection
$G_h$ can be extended to a smooth function on $B\left(0,
\frac{1}{2}\right)\bs\{\pm he_{n+1}\}$ with $G_h(\bar x, x_{n+1})=-G(\bar x,
-x_{n+1})$ for $x_{n+1}>0$. For $j$ large let $G^j_h(X)=G_n(X+\s_je_{n+1})$
be defined on
$B\left(\frac{1}{2}, -\s_je_{n+1}\right)\bs \{(\s_j\pm h)e_{n+1}\}$. We
denote by $B_{1/2}=B\left(0, \frac{1}{2}\right)$ and by
$B^j_{1/2}=B\left(\frac{1}{2}; -\s_je_{n+1}\right)$. We may assume that
$\H^n(\p B^j_{1/2}\cap\p\{G_j>0\})=0$. Green's formula ensures that
\begin{equation}\label{eqn3.64b}-\int_{B^j_{1/2}}\langle\nabla G_j, \nabla
G^j_h\rangle = \int_{\p B^j_{1/2}}G_j\p_\nu G^j_h-G_j(-(h+\s_j)e_{n+1}),
\end{equation}
where $\p_\nu G^j_h=\langle\nabla G^j_h,\nu\rangle$, and $\nu$ denotes the
inward pointing unit normal to $\p B'_{1/2}$. On the other hand
\begin{equation}\label{eqn3.65b}
-\int_{\p B^j_{1/2}}\langle\nabla G_j,\nabla G^j_h\rangle =
\int_{\p\{G_j>0\}\cap B^j_{1/2}} h_j G^j_h d\H^n.
\end{equation}

Let $\nu_j$ denote the inward point unit normal to $\po_j=\p\{G_j>0\}$ then
by Green's formula we have
\begin{equation}\label{eqn3.66b}
\int_{B^j_{1/2}\cap\p\{G_j>0\}}\langle G^j_he_{n+1}-x_{n+1}\nabla G^j_h,
\nu_j\rangle  d\H^n= (\s_j+h) + \int_{B^j_{1/2}\cap\{G_j>0\}}x_{n+1}\p_\nu G^j_h.
\end{equation}
Combining (\ref{eqn3.64b}), (\ref{eqn3.65b}) and (\ref{eqn3.66b}) we obtain
\begin{eqnarray}\label{eqn3.67b}
&&\int_{B^j_{1/2}\cap\p\{G_j>0\}}x_{n+1}\p\nu_{j} G^j_h d\H^n  \\
& = &
\int_{B_{1/2}\cap\p\{G_j>0\}}(h_j+\langle e_{n+1},
\nu_j\rangle)G^j_h  d\H^n\nonumber \\
&& -
\int_{\p B^j_{1/2}\cap\{G_j>0\}}(x_{n+1}+G_j)\p_\nu G^j_h +
G_j(-(h+\s_j)e_{n+1})-(\s_j+h)\nonumber\\
&= & \int_{B_{1/2}\cap \p\{G_j>0\}}\left(\frac{h_j}{1-\tau_j}+\langle
e_{n+1},\nu_j\rangle\right)G^j_h  d\H^n \nonumber \\
&& - \tau_j\int_{B_{1/2}\cap\p\{G_j>0\}} h_jG^j_h d\H^n+
G_j(-(h+\s_j)e_{n+1})-(\s_j+h) \nonumber \\
&& - \int_{\p B^j_{1/2}\cap\{G_j>0\}}(x_{n+1}+G_j)\p_\nu G^j_n \nonumber \\
& = & \int_{B_{1/2}\cap\p\{G_j>0\}} \left(\frac{h_j}{1-\tau_j} + \langle
e_{n+1}, \nu_j\rangle\right)G^j_h  d\H^n
+ (1+\tau_j)G_j(-(h+\s_j)e_{n+1})-(\s_j+h) 
\nonumber \\
&& - \int_{\p B^j_{1/2}\cap\{G_j>0\}}(x_{n+1}+G_j(1+\tau_j))\p_\nu G^j_h.
\nonumber
\end{eqnarray}
Since $\s_j-h<-\s_j$ and $G_j\in F(\s_j, \s_j; \tau_j)$ in $B(0,1)$ in
direction $e_{n+1}$, then $G^j_h\le 0$ on $\p\{G_j>0\}\cap B^j_{1/2}$.
Furthermore since $h_j(0)=1$, by (\ref{eqn3.54}) $h_j\ge 1-\tau_j$ on
$B^j_{1/2}\cap \p\{G_j>0\}$.

Thus
\begin{equation}\label{eqn3.66t}
\int_{B_{1/2}\cap\p\{G_j>0\}}\left(\frac{h_j}{1-\tau_j} + \langle
e_{n+1},\nu_j\rangle\right) G^j_h\le 0.
\end{equation}
Since $G_j(0)=0$ (\ref{eqn3.54}) ensures that
\begin{equation}\label{eqn3.67t}
|G_j(-(h+\s_j)e_{n+1})| \le  \sup_{B(0,1)}|\nabla G_j|(h+\s_j) 
 \le  (1+\tau_j)(h+\s_j).
\end{equation}
Hence
\begin{equation}\label{eqn3.68t}
(1+\tau_j)G_j(-(h+\s_j)e_{n+1}) - (\s_j+h) \le 3\tau_j(h+\s_j)
\end{equation}
Since $\{G_j>0\}\subset\{x_{n+1}<\s_j\}$, by (\ref{eqn3.54}) for
$x_{n+1}\le\s_j$ we have in $B(0,1)$
\begin{equation}\label{eqn3.69t}
G_j(\bar x, x_{n+1}) = |G_j(\bar x, x_{n+1}) - G_j(\bar x, \s_j)|\le
(1+\tau_j)(\s_j-x_{n+1})
\end{equation}
which yields
\begin{equation}\label{eqn3.70t}
x_{n+1}\le x_{n+1} + G_j(1+\tau_j)\le (1-(1+\tau_j)^2)x_{n+1} +
(1+\tau_j)^2\s_j.
\end{equation}

Thus
\begin{equation}\label{eqn3.71t}
0 \le x_{n+1} + (1+\tau_j)G_j\le (1+\tau_j)^2\s_j\mbox{ for }x_{n+1}\in
[0,\s_j]
\end{equation}
\begin{equation}\label{eqn3.72t}
-\s_j\le x_{n+1} + (1+\tau_j)G_j\le (1+\tau_j)\s_j\mbox{ for }x_{n+1}\in
[-\s_j,0].
\end{equation}
Since $G_j\in F(\s_j, \s_j;\tau_j)$ in $B(0,1)$ in direction $e_{n+1}$ with
$h_j(0)=1$ then
\begin{eqnarray}\label{eqn3.73t}
x_{n+1}+G_j(1+\tau_j) & \ge & x_{n+1}+(1+\tau_j)(-x_{n+1}-\s_j) \\
& \ge & -\tau_jx_{n+1} - \s_j(1+\tau_j)\ge -\s_j(1+\tau_j)\mbox{ for
}x_{n+1} \le -\s_j \nonumber
\end{eqnarray}

We combine the fact that $\p_\nu G^j_h\ge 0$ with (\ref{eqn3.71t}),
(\ref{eqn3.72t}) and (\ref{eqn3.73t}) and obtain that
\begin{equation}\label{eqn3.74t}
-\int_{\p B^j_{1/2}\cap\{G_j>0\}} (x_{n+1} + (1+\tau_j)G_j)\p_\nu G^j_h\le
\s_j(1+\tau_j) \int_{\p B^j_{1/2}\cap\{G_j>0\}\cap\{x_{n+1}<0\}} \p_\nu
G^j_n.
\end{equation}
Combining (\ref{eqn3.67b}), (\ref{eqn3.66t}), (\ref{eqn3.68t}),
(\ref{eqn3.74t}),  the fact that $\s^{-2}_j\tau_j\le 1$ for $j$ large
enough, and that $1\ge h>2\s_j$ we conclude that
\begin{equation}\label{eqn3.75t}
\frac{1}{\s_j} \int_{B^j_{1/2}\cap\p\{G_j>0\}} x_{n+1}\p_{\nu_j}G^j_h \le
9\s_j + 2 \int_{\p B^{1/2}_{1/2}\cap \{G_j>0\} \cap \{x_{n+1}<0\}} \p_\nu
G^j_n.
\end{equation}

Thus
\begin{equation}\label{eqn3.76t}
\limsup_{j\to\infty} \frac{1}{\s_j} \int_{B^j_{1/2}\cap\p\{G_j>0\}}
x_{n+1}\p_{\nu_j}G^j_h \le 2\int_{\p B_{1/2}\cap\{x_{n+1}\le 0\}} \p_\nu
G_h\le Ch.
\end{equation}
The rest of the argument is identical to the one that appears in \cite{KT1}
in the proof of Lemma 0.9.
\qed
\end{proof}

\begin{lemma}[Lemma 7.7 and Lemma 7.8 \cite{AC}]\label{lem3.7}
The function $f$ introduced in Lemma \ref{lem3.4} is Lipschitz in
$B'_{1/16}$ with Lipschitz constant that only depends on $n$. Furthermore
there exists a large constant $C=C(n)>0$ such that for
any given $\theta\in(0,1)$ there exists $\eta=\eta(\theta)>0$ and
$l\in\RR^n\times\{0\}$ with $|l|\le c$ so that
\begin{equation}\label{eqn3.77t}
f(y)\le \langle l,y\rangle + \frac{\theta}{2}\eta\mbox{ for }y\in B'_\eta.
\end{equation}
\end{lemma}

The proof of this lemma basically appears in \cite{AC} and \cite{KT1}

Now we indicate how the last 2 lemmata yield a contradiction in the proof of
Lemma \ref{lem3.2}. Recall that by assuming that the statement in Lemma
\ref{lem3.2} is false we can construct sequences of function $\{G_j\}$ and
$\{h_j\}$ satisfying (\ref{eqn3.53}), (\ref{eqn3.54}), (\ref{eqn3.55}) and
(\ref{eqn3.56}). From them as in (\ref{eqn3.57}), (\ref{eqn3.58}) and
Lemmas \ref{lem3.4}, \ref{lem3.5}, \ref{lem3.6} and \ref{lem3.7} we can
produce a subharmonic Lipschitz function $f$ on $B'_{1/16}$ satisfying
(\ref{eqn3.77t}). Recall that by Lemma \ref{lem3.4} and Corollary
\ref{cor3.1} $f$ is uniform limit of the functions $f^+_j$ defined in
(\ref{eqn3.57}). Therefore Lemma \ref{lem3.7} yields that for
$\theta\in(0,1)$ there exists $\eta>0$ so that for $j$ large enough
\begin{equation}\label{eqn3.78t}
f^+_j(y)\le\langle l, y\rangle +\theta\eta\mbox{ for }y\in B'_\eta,
\end{equation}
which by definition means that
\begin{equation}\label{eqn3.79t}
G_j(X)=0\mbox{ for }X=(\bar x, x_{n+1}) \in B(0,\eta)\mbox{ with
}x_{n+1}>\s_j\langle l,\bar x\rangle+\theta\eta\s_j.
\end{equation}

Let $\bar\nu=(1+\s^2_j|l|^2)^{-1/2}(-\s_jl,1)$ (\ref{eqn3.79t}) implies that
\begin{equation}\label{eqn3.80t}
G_j(X)=0\mbox{ for }X\in B(0,\eta) \mbox{ with }\langle X, \bar\nu\rangle
\ge \frac{\theta\eta\s_j}{(1+\s^2_j|l|^2)^{1/2}}\ge 2\theta\eta\s_j
\end{equation}
for $j$ large enough. But (\ref{eqn3.54}) and (\ref{eqn3.80t}) state that
$G_j\in F(2\theta\eta_j, 1; \tau_j)$ in $B(0,\eta)$ in direction $\bar\nu$.
This
 contradicts statement (\ref{eqn3.56}) in the case that
$\theta=\frac{\theta_0}{2}$, which concludes the proof of Lemma \ref{lem3.2}
and thus that of the Theorem \ref{thm1}.

%% file: stab-apps.tex
\section{Applications}

\begin{lemma}\label{lem4.1}
Assume that $\O\subset\RR^{n+1}$ satisfies (\ref{eqn1}). Then there exist
$\e_0>0$ and $r_0>0$ such that if
\begin{equation}\label{eqn4.1}
\sup_{\po} |\log h|<\e_0
\end{equation}
then for $Q\in\po$ and $r\in (0,r_0)$
\begin{equation}\label{eqn4.2}
C^{-1}_nr^n\le\H^n(\po\cap B(Q,r))\le C_nr^n,
\end{equation}
where $C_n$ is a constant that only depends on $n$, i.e. $\po$ is Ahlfors
regular.
\end{lemma}

\begin{proof}
Let $\s\in\left(0,\frac{1}{4}\right)$ be small enough in Theorem
\ref{thm1} then there exists $\e_1>0$ such that if $\sup_{\po}|\log h|<\e_1$,
then $\po$ is $\s$-Reifenberg flat. This ensures that there exists
$\rho_1>0$ so that for $Q\in\po$ and $r<\rho_1$
\begin{equation}\label{eqn4.3}
\H^n(\po\cap B(Q,r))\ge (1+\s)^{-1}\o_nr^n\ge \frac{1}{2}\o_nr^n
\end{equation}
(for the proof see Remark 2.2 in \cite{KT2}). By Lemma \ref{lem3.3} 
there exists
$0<\e_2<\e_1$ so that if $\sup_{\po}|\log h|<\e$ with $0<\e<\e_2$ there
exists $\rho_\e=\rho>0$ such that for $Q\in\po$,\newline
$G\in F(\s,\s;(e^{2\e}-1)^{1/4})$ in $B(Q, \rho_\e)$. Thus in particular for
$r<\min\{\rho_\e, \rho_1\}$
\begin{equation}\label{eqn4.4}
\sup_{B(Q,r)}|\nabla G|  \le  h(Q)(1+(e^{2\e}-1)^{1/4}) 
\le  e^\e(1+(e^{2\e}-1)^{1/4}). 
\end{equation}
Hence
\begin{eqnarray}\label{eqn4.5}
\H^n(\po\cap B(Q,r)) & = & \int_{B(Q,r)\cap\po} h
\frac{1}{h}d\H^n \\
& \le & e^{+\e}\int_{B(Q,r)\cap\po} hd\H^n\le e^{+\e}\int_{\po}\varphi
hd\H^n \nonumber \\
& \le & -e^{+\e}\int_\O\langle\nabla \varphi,\nabla G\rangle d\H^{n+1}, \nonumber
\end{eqnarray}
for any non-negative $\varphi\in C^\infty_c(\RR^{n+1})$ such that 
$\varphi\equiv 1$ on $B(Q,r)$ and $0\not\in{\rm support}\,\varphi$.

In particular if $\varphi$ is chosen so that $\varphi\in
C^\infty_c(B(Q,2r))$ for $r<\frac{1}{2}\min\{\rho_\e, \rho_1\}$ and
$|\nabla\varphi|<2/r$, (\ref{eqn4.4}) and (\ref{eqn4.5}) yield for
$\e>0$ small enough
\begin{equation}\label{eqn4.6}
\H^n(B(Q,r)\cap\po)  \le 
e^{\e}\frac{2}{r}e^{\e}(1+(e^{2\e}-1)^{1/4})\o_{n+1}r^{n+1} 
\le  4\o_{n+1}r^n. 
\end{equation}
Choosing $\e_0=\min\left\{\frac{\e_2}{2}, \frac{1}{4}\right\}$ and
$\rho_0=\frac{1}{2}\min\{\rho_{\e_0}, \rho_1\}$ we conclude that
(\ref{eqn4.2}) holds.
\end{proof}

\begin{corollary}\label{cor4.1}
Assume that $\O\subset\RR^{n+1}$ satisfies (\ref{eqn1}). Then given $\d>0$
small enough there exists $\e>0$ such that if
\begin{equation}\label{eqn4.7}
\sup_{\po}|\log h|<\e
\end{equation}
then $\O$ is a $\d$-Reifenberg flat chord arc domain.
\end{corollary}

\begin{proof}
By Theorem \ref{thm1}, $\po$ is $\d$-Reifenberg flat provided $\e>0$ is
small enough. Since $\O$ is bounded and $B(0,r_1)\subset\O\subset B(0, R_2)$ 
it is easy to show that it satisfies the
separation property. Therefore $\O$ is a $\d$-Reifenberg flat domain and
for $\d>0$ small enough it is also NTA (see \cite{KT2}). Moreover if
$\e<\e_0$ Lemma \ref{lem4.1} ensures that for $r\in(0,r_0)$
(\ref{eqn4.2}) holds. Since $\O$ is bounded it is easy to see that for $r\in
(0, \diam \O)$, (\ref{eqn4.2}) also holds with a constant that only depends
on $n$, and $\frac{\diam\O}{\rho_0}$. Thus $\O$ is a chord arc domain.
\qed
\end{proof}

The crucial information contained in Lemma \ref{lem4.1} and Corollary
\ref{cor4.1} is that bounded domains which are sets of locally finite
perimeter and satisfy (\ref{eqn1}) belong to a family of chord arc domains
with uniform constants.

\begin{corollary}\label{cor4.2}
Assume that $\O\subset\RR^{n+1}$ satisfies (\ref{eqn1}). There exists
$\e_1>0$ so that if $\sup_{\po}|\log h|<\e_1$ and $\log h\in \vmo(\po)$
(resp.\ $\log h\in C^{k,\a}(\po)$) then $\O$ is a chord arc domain with
vanishing constant (resp.\ $\O$ is a $C^{k+1,\a}$ domain).
\end{corollary}

\begin{proof}
By choosing $\e_1>0$ small enough Corollary \ref{cor4.1} ensures that $\O$
is a $\d$-Reifenberg flat chord arc domain. Choosing $\d>0$ as in the
statement of the Main Theorem in \cite{KT3} we conclude that if $\log
h\in\vmo$ then $\vecn \in\vmo(\po)$. Choosing $\d>0$ as in the statement
of Alt and Caffarelli's theorem we conclude that if $\log h\in C^{k,\a}$
then $\O$ is a $C^{k+1,\a}$ domain.
\qed
\end{proof}

\begin{corollary}\label{cor4.2A}
Assume that $\O\subset\RR^{n+1}$ satisfies (\ref{eqn1}). 
There exists $\e_2>0$ so
that it $\sup_{\po}|\log h|<\e_2$ and $\log h\in C^{0,\a}$ there exists
a homeomorphism
$\psi:B(0,R_1)\rightarrow\O$ where $\psi$ and $\psi^{-1}$ are $C^{1,\a}$.
\end{corollary}

\begin{proof}
By the work in \cite{AC} and Corollary \ref{cor4.1} we know that there
exists $\d>0$ and $\e>0$ depending on $\d>0$ so that if $\sup_{\po}|\log
h|<\e$ and $\log h\in C^{0,\a}$ then $\O$ is a $C^{1,\a}$ domain. Moreover
using the proof of Theorem 8.1 in \cite{AC} and (\ref{eqn3.6}) above we
conclude that
\begin{equation}\label{eqn4.1A}
\left|\vecn(Q)-\frac{Q}{|Q|}\right|<\d.
\end{equation}

Here $\vecn(Q)$ denotes the outward unit normal to $\po$. Since $\O$ is a
bounded $C^{1,\a}$ domain there exists 
$r\in\left(0, \frac{R_1}{8}\right)$ so that for
$Q\in\po\cap B(Q,r)$ can be written as the area below the graph of a
$C^{1,\a}$ function (with small $C^{1,\a}$ norm 1 over the $n$-plane through
$Q$ and orthogonal to $\vecn(Q)$. Inequality (\ref{eqn4.1A}) guarantees that
$\O\cap B(Q,r)$ can also be seen as the area below the graph of a $C^{1,\a}$
function (with $C^{1,\a}$ norm less than $C\d$) over the $n$-plane through
$Q$ and orthogonal to $\frac{Q}{|Q|}$. This implies that the spherical
projection $S:\po\to B(0,R_1)$ $S(Q)=R_1\frac{Q}{|Q|}$ is a 1-1 map.
Moreover since $B(Q,R_1)\subset \O$, $S$ is onto and Lipschitz on $\po$. In
particular $\O$ is star shaped with respect to the origin.

Since $S$ is smooth on $\RR^{n+1}\bs B\left(0, \frac{R_1}{4}\right)$ and
$\po$ is a $C^{1,\a}$ submanifold it is clear that $S$ is a $C^{1,\a}$ map
from $\po$ onto $B(0,R_1)$, and $S^{-1}$ is a $C^{1,\a}$ map from $\p
B(0,R_1)$ onto $\po$. For $X\in \O\bs B\left(0, \frac{R_1}{4}\right)$ there
exists a unique $Q_X\in \po$ so that $\frac{X}{|X|} = \frac{Q_X}{|Q_X|}$.
The previous remark ensures that the map that to $X\in\O\bs B\left(0,
\frac{R_1}{4}\right)$ associates $Q_X$ is a $C^{1,\a}$ map. Our goal is to
construct a homeomorphism $\Phi:\O\rightarrow B(0,R_1)$, such that $\Phi$ and
$\Phi^{-1}$ are $C^{1,\a}$. Let $X\in\O$ and define
\begin{equation}\label{eqn4.2A}
g(t) = \left\{\begin{array}{ll}
t & t\in\left[0,\frac{R_1}{4}\right] \\
\frac{R_1-|Q_X|}{\left(|Q_X|-\frac{R_1}{4}\right)^2}\left(t-\frac{R_1}{4}\right)^2 + t &
\mbox{for }t\in\left[\frac{R_1}{4}, |Q_X|\right]. \end{array}\right.
\end{equation}
In particular $g\in C^{1,1}([0,|Q_X|])$, $g(0)=0$ and $g(|Q_X|)=R_1$. Moreover
since $|Q_X|\ge R_1$, for $\e<\frac{1}{64}$, $g'>0$ on $[0, |Q_X|]$ thus $g$
is $1-1$ and maps $[0, |Q_X|]$ onto$ [0, R_1]$. For $X\in\O$
define
\begin{equation}\label{eqn4.3A}
\Phi(X) = g(|X|)\frac{X}{|X|} - \left\{ \begin{array}{lcl}
X & \mbox{for} & X\in B\left(0, \frac{R_1}{4}\right) \\
\left(\frac{R_1-|Q_X|}{\left(|Q_X|-\frac{R_1}{4}\right)^2}\left(|X| -
\frac{R_1}{4}\right)^2 + |X|\right) \frac{X}{|X|} & \mbox{for} & X\in\O\bs
B\left(0, \frac{R_1}{4}\right). \end{array}\right.
\end{equation}
Note that $\Phi$ is a $C^{1,\a}$ map.
For $Y\in B(0,R_1)\subset\O$ there exists a unique $Q_Y\in\po$. Since $g$ is
a bijection there exists a unique $t\in[0, |Q_Y|]$ so that $|Y|=g(t)$. Since
$\O$ is star-shaped with respect to the origin there exists $X\in\O$, such
that $X=t\frac{Q_Y}{|Q_Y|}$. This implies that $\Phi(X)=Y$. If
$\Phi(X)=\Phi(X')\Rightarrow g(|X|)=g(|X'|)$ and $\frac{X}{|X|} =
\frac{X'}{|X'|}$. Since $g$ is $1-1$, $|X|=|X'|$ which yields  $X=X'$. Thus
$\Phi:\O\to B(0, R_1)$ is a $C^{1,\a}$ bijection. It is easy to check that
$\Phi^{-1}$ is also $C^{1,\a}$.
\qed
\end{proof}

\begin{lemma}\label{lem4.2}
Assume that $\O\subset\RR^{n+1}$ satisfies (\ref{eqn1}). Given $\d>0$ there
exists $\e>0$ such that if $\sup_{\po}|\log h|<\e$ with $\e<\e_0$ then there
exists $\rho_\e>0$ such that for $r\in (0,\rho_\e)$ and $Q\in\po$
\begin{equation}\label{eqn4.8}
\frac{H^n(B(Q,r)\cap\po)}{\o_nr^n}\le (1+\d).
\end{equation}
\end{lemma}

\begin{proof}
Let $\varphi\in C^\infty_c(\RR^{n+1})$ such that $0\not\in\supp\varphi$ then
\begin{equation}\label{eqn4.9}
-\int_\O\langle\nabla\varphi,\nabla G\rangle = \int_{\po}\varphi hd\s.
\end{equation}
By choosing $\varphi$ as an approximation of $\chi_{B(Q,r)}$ we obtain after
passing to the limit that for a.e. $r>0$ with $r<\frac{R_1}{4}$
\begin{equation}\label{eqn4.10}
\int_{\po\cap B(Q,r)}hd\H^n = \int_{\p B(Q,r)\cap\O} \left\langle \nabla G,
\frac{X-Q}{|X-Q|}\right\rangle d\H^n.
\end{equation}
For the details of this computation see \cite{KT4} section 3.

Let $\d'=\d'(\d)\in (0,1)$ and choose $\e'_0\in\left(0,\frac{1}{4}\right)$
so that if $\sup|\log h|<\e'$ for $\e'\in(0,\e'_0)$ then $\O$ is a
$\d'$-Reifenberg flat chord arc domain (see Corollary \ref{cor4.1}) and
$G\in F\left(\frac{\d'}{2}, \frac{\d'}{2}, (e^{2\e'-1})^{1/4}\right)$ in
$B\left(Q, \left(\frac{\eta}{2}\right)^k\rho'\right)$ for all $Q\in\po$,
$k\ge 1$ where $\rho'=\sqrt{2}\sqrt{e^{2\e'-1}}R_1$ and $\eta\in\left(0,
\frac{1}{4}\right)$ (see Lemma \ref{lem3.3} and the proof of Theorem
\ref{thm1}, namely (\ref{eqn3.29}) and (\ref{eqn3.30})).

By Lemma \ref{lem2}, $B(0, R_1)\subset\O\subset B(0,R_2)$ with
$e^{-\e'}\le\s_nR^n_1\le\s_nR^n_2\le e^{\e'}$. Note that since $G\in
F\left(\frac{\d'}{2}, \frac{\d'}{2}, (e^{2\e'}-1)^{1/4}\right)$ in
$B\left(Q, \left(\frac{\eta}{2}\right)^k\rho'\right)$ for $k\ge 1$ then
$G\in F(\d', \d', (e^{2\e'}-1)^{1/4})$ in $B(Q,r)$ for $r\in (0, \rho')$.
Thus there exists $\vecnQr\in \SS^n$ so that
\begin{equation}\label{eqn4.11}
G(X)=0\qquad\mbox{for }\qquad\langle X-Q; \vecnQr\rangle\le-\d'r
\end{equation}
and
\begin{equation}\label{eqn4.12}
G(X)\ge h(Q)\left[\langle X-Q; \vecnQr\rangle-\d'r\right] \ \mbox{ for
}\ \langle X-Q; \vecnQr\rangle\ge\d'r
.
\end{equation}
To estimate the term in the right hand side of (\ref{eqn4.10}) consider
\begin{eqnarray}\label{eqn4.13}
0 & \le & \int_{\p B(Q,r)\cap\O}\left\langle\nabla G;
\frac{X-Q}{|X-Q|}\right\rangle d\H^n\\
& \le & \int_{\p B(Q,r)\cap\{x+t\vecnQr: t\ge 2\sqrt{\d'}r\}}
\left\langle\nabla G(X); \frac{X-Q}{|X-Q|}\right\rangle d\H^n \nonumber \\
&&\qquad + \int_{\p
B(Q,r)\cap\{x+t\vecnQr: -\d'r\le t\le 2\sqrt{\d'}r\}}|\nabla G|d\H^n
\nonumber 
\end{eqnarray}
Here the decomposition $x+t\vecnQr$ means that $x\in L(Q,r)$ where
$L(Q,r)$ is an $n$-plane through $Q$, orthogonal to $\vecnQr$.

Given our choice of $r$, Lemma \ref{lem3.3} guarantees that
\begin{equation}\label{eqn4.14}
\sup_{B(Q,r)}|\nabla G| \le  h(Q) (1+(e^{2\e'1}-1)^{1/4}) 
 \le  h(Q)(1+2(\e')^{1/4}), 
\end{equation}
for $\e'_0$ small enough.

Using (\ref{eqn4.14}) a simple computation yields
\begin{equation}\label{eqn4.15}
\int_{\p B(Q,r)\cap\{x+t\vecnQr; -\d'_r\le t\le 2\sqrt{\d'}r\}}
|\nabla G|d\H^n \le C_n\sqrt{\d'}r^n.
\end{equation}
Combining (\ref{eqn4.11}), (\ref{eqn4.12}) and (\ref{eqn4.14}) we have for
$X\in B(Q,r)$, $X=x+t\vecnQr$ with $t\ge 2\d'r\ge 2\sqrt{\d'}r$
\begin{equation}\label{eqn4.16}
h(Q)(t-\d'r)\le G(X)\le h(Q)(1+2(\e')^{1/4})(t+\d'r).
\end{equation}
Note that for such $X$, if $d(X)$ denotes the distance from $X$ to $\po$
then
\begin{equation}\label{eqn4.17}
r\ge d(X)\ge t-\d'r\ge\frac{t}{2}.
\end{equation}

As in (\ref{eqn3.23}) and (\ref{eqn3.24}) we have that
\begin{eqnarray}\label{eqn4.18}
\nabla G(X) & = & \frac{-2^{n+2}}{\o_{n+1}d(X)^{n+2}} \int_{\p B\left(X,
\frac{d(X)}{2}\right)} G(\zeta)(X-\zeta)d\zeta \\
& = & -\frac{2^{n+2}}{\o_{n+1}d(X)^{n+2}} \int_{\p B\left(X,
\frac{d(X)}{2}\right)} (G(\zeta) - h(Q)\widetilde t_\zeta)(X-\zeta)d\zeta
\nonumber \\
& - & \frac{2^{n+2}}{\o_{n+1}d(X)^{n+2}} \int_{\p B\left(X,
\frac{d(X)}{2}\right)} h(Q)\widetilde t_\zeta(X-\zeta)d\zeta, \nonumber
\end{eqnarray}
where $\widetilde t_\zeta = \langle\zeta-Q, \overrightarrow{n}_{Q,2r}\rangle$. 
Note that
if $\zeta\in\p B\left(X, \frac{d(X)}{2}\right)$, then $\zeta\in B(Q,2r)$.
The first equality in (\ref{eqn4.18}) applied to the function 
$\widetilde t_\zeta$ guarantees that
\begin{equation}\label{eqn4.19}
\frac{2^{n+2}}{\o_{n+1}d(X)^{n+2}} \int_{\p B\left(X, \frac{d(X)}{2}\right)}
h(Q)\widetilde t_\zeta(X-\zeta)d\zeta = h(Q)\overrightarrow{n}_{Q,2r}.
\end{equation}

Since $|\widetilde t_\zeta|\le 2r$ using (\ref{eqn4.16}) we have that
\begin{eqnarray}\label{eqn4.20}
|\nabla G(X)-h(Q)\overrightarrow{n}_{Q,2r}| 
& \le & \frac{C_n}{d(X)^{n+1}} \int_{\p
B\left(X, \frac{d(X)}{2}\right)} |G(\zeta) - h(Q)\widetilde t_\zeta|d\zeta \\
& \le & \frac{C_nh(Q)}{d(X)^{n+1}} \int_{\p B\left(X, \frac{d(X)}{2}\right)}
((\e')^{1/4}(|\widetilde t_\zeta|+\d'r)+\d'r)d\zeta \nonumber \\
& \le & \frac{C_nh(Q)}{d(X)} ((\e')^{1/4}r+\d'r) \le C_n\frac{h(Q)}{t}
((\e')^{1/4}r + \d'r). \nonumber
\end{eqnarray}

Using (\ref{eqn4.17}) and (\ref{eqn4.20}) we can estimate the remaining term
in (\ref{eqn4.13}). Namely
\begin{eqnarray}\label{eqn4.21}
&&\kern1in\int_{\p B(Q,r)\cap\{x+t\vecn_{Q,r}: t\ge 2\sqrt{\d'}r\}} 
\left\langle
\nabla G(X); \frac{X-Q}{|X-Q|}\right\rangle\H^n \\
& \le &  h(Q) \int_{\p
B(Q,r)\cap \{x+t\vecn_{Q,r}: t\ge 2\sqrt{\d'}r\}}
\left\langle\vecn_{Q,2r}, \frac{X-Q}{|X-Q|}\right\rangle d\H^n + C_nh(Q)
\frac{(\e')^{1/4}+\d'}{\sqrt{\d'}} r^n.\nonumber
\end{eqnarray}

Choosing $\e'_0>0$ so that $\e'_0\le(\d')^4$, and recalling that $h(Q)\le
e^{\e'}\le 2$ (\ref{eqn4.21}) becomes
\begin{eqnarray}\label{eqn4.22}
&&\kern1in\int_{\p B(Q,r)\cap\{x+t\vecnQr: t\ge 2\sqrt{\d'}r\}}\left\langle 
\nabla G; \frac{X-Q}{|X-Q|}\right\rangle d\H^n \\
& \le & h(Q) \int_{\p B(Q,r)\cap \{x +
\widetilde t\vecn_{Q,2r}: \widetilde t\ge 0\}} \left\langle \vecn_{Q, 2r},
\frac{X-Q}{|X-Q|}\right\rangle d\H^n \nonumber \\
&&\quad + 2\H^n\left(\p B(Q,r)\cap\left(\{x+\widetilde t\vecn_{Q,2r}: 
\widetilde
t\ge 0\}\Delta \{x+t\vecn_{Q,r}:t\ge 2\sqrt{\d'}r\}\right)\right) +
C_n\sqrt{\d'}r^n.\nonumber
\end{eqnarray}

A simple computation shows that the angle between $\vecn_{Q,2r}$ and
$\vecn_{Q,r}$ is less than $C\d'$. This fact combined with (\ref{eqn4.10})
applied to the function $\widetilde t\vecn_{Q,2r}$ instead of $G$ and
(\ref{eqn4.22}) implies
\begin{eqnarray}\label{eqn4.23}
&&\kern-1in\int_{\p B(Q,r)\cap\{x+t\vecn_{Q,r}, t\ge 2\sqrt{\d'}r\}}\left\langle \nabla
G, \frac{X-Q}{|X-Q|}\right\rangle d\H^n \\
& \le & h(Q) \int_{L(Q,2r)\cap B(Q,r)}d\H^n + C_n\sqrt{\d'}r^n \nonumber
\end{eqnarray}

Combining (\ref{eqn4.10}), (\ref{eqn4.13}), (\ref{eqn4.15}) and
(\ref{eqn4.23}) plus the fact that $e^{-\e'}\le h(P)\le e^{\e'}$ for
$P\in\po$ we conclude for $r\le \sqrt{2}\sqrt{e^{2\e'}-1}R_1$
\begin{eqnarray}\label{eqn4.24}
\H^n(B(Q,r)\cap\po) & = & \int_{B(Q,r)\cap\po} h(P) h^{-1}(P)d\H^n \\
& \le & e^{\e'}\int_{B(Q,r)\cap\po}hd\H^n \le
e^{2\e'}\o_nr^n+C_n\sqrt{\d'}r^n. \nonumber
\end{eqnarray}

By our choice of $\e'_0>0$ (so that $\e'_0\le (\d')^4$) we have that for
$r\le\sqrt{2}\sqrt{e^{2\e'}-1}R_1$
\begin{equation}\label{eqn4.25}
\H^n(B(Q,r)\cap\po)\le \o_nr^n(1+C_n\sqrt{\d'}).
\end{equation}
Choosing $\d'>0$ so that $C_n(\d')^{1/2}=\d$, and $\e_0$ the corresponding
$\e'_0$ we have proved the statement of Lemma \ref{lem4.2}.
\qed
\end{proof}

\begin{corollary}\label{cor4.3}
Assume that $\O\subset\RR^{n+1}$ satisfies (\ref{eqn1}). Given $\d>0$ there
exists $\e>0$ such that if $\sup_{\po}|\log h|<\e$ then $\O$ is
a $\d$-chord arc domain.
\end{corollary}

\begin{proof}
From the proof of Lemma \ref{lem4.1} (see (\ref{eqn4.3})) and Lemma
\ref{lem4.2} we have that given $\d>0$ there exist $\e>0$ and $\rho>0$ so
that if $\sup_{\po}|\log h|<\e$ then for $r\in(0,\rho)$ and $Q\in\po$
\begin{equation}\label{eqn4.26}
(1+\d)^{-1}\le \frac{H^n(\po\cap B(Q,r))}{\o_nr^n}\le 1+\d.
\end{equation}
By Theorem \ref{thm1} we also know that $\rho>0$ can be chosen so that
\begin{equation}\label{eqn4.27}
\theta(Q,\rho)\le\d.
\end{equation}
\qed
\end{proof}

This is a straightforward consequence of Corollary 4.5 (for the
proof see \cite[\S2]{KT2}).

\begin{corollary}\label{cor4.4}
Assume that $\O\subset\RR^{n+1}$ satisfies (\ref{eqn1}). Given $\d>0$ there
exist $\e>0$ and $\rho>0$ such that if $\sup_{\po} |\log h|<\e$ then 
\begin{equation}\label{eqn4.28}
\|\vecn\|_\ast(\rho) = \sup_{Q\in\po} \sup_{0<r<\rho}
\left(\notint_{B(Q,r)\cap\po} |\vecn-\vecn_{Q,r}|^2d\s\right)^{1/2} \le \d.
\end{equation}
\end{corollary}

%% file: stab-refs.tex
\vskip 1in

D. Preiss: Department of Mathematics, University College London, \\
dp@math.ucl.ac.uk\\

T. Toro: Department of Mathematics, University of Washington,\\
toro@math.washington.edu